\RequirePackage{fix-cm}
\documentclass[12pt,a4paper,onecolumn,twoside]{article}

\usepackage[utf8]{inputenc}
\usepackage{amsmath,amssymb,amsfonts,amsthm,mathrsfs}
\usepackage{mathtools}
\usepackage{graphicx}
\usepackage{tikz}
\usetikzlibrary{calc,arrows.meta}
\definecolor{ink}{HTML}{293A47}
\definecolor{accent}{HTML}{166581}
\definecolor{guide}{HTML}{70828C}
\definecolor{wash}{HTML}{F3F6F7}
\definecolor{imagefill}{HTML}{C3E1E7}
\usepackage{geometry}
\usepackage{fancyhdr}
\usepackage{titlesec}
\usepackage{setspace}
\usepackage{indentfirst}
\usepackage[runin]{abstract}

\abslabeldelim{.\hspace{0.35em}}
\AtBeginDocument{\setlength{\absparindent}{0pt}}
\usepackage[numbers,sort]{natbib}
\usepackage{enumitem}
\usepackage[hidelinks,hyperfootnotes=false]{hyperref}

\setenumerate[1]{label=\textnormal{(\arabic*)},itemsep=-4pt,partopsep=0pt,parsep=\parskip,topsep=5pt}
\setitemize[1]{itemsep=-4pt,partopsep=0pt,parsep=\parskip,topsep=5pt}
\setdescription{itemsep=-4pt,partopsep=0pt,parsep=\parskip,topsep=5pt}
\allowdisplaybreaks
\titleformat{\section}[block]{\normalfont\fontsize{14pt}{17pt}\selectfont\bfseries\raggedright}{\thesection}{1em}{}
\titleformat{\subsection}[block]{\normalfont\fontsize{12pt}{14.5pt}\selectfont\bfseries\raggedright}{\thesubsection}{1em}{}
\titleformat{\subsubsection}[block]{\normalfont\fontsize{12pt}{14.5pt}\selectfont\bfseries\raggedright}{\thesubsubsection}{1em}{}
\titleformat{\paragraph}[block]{\normalfont\fontsize{12pt}{14.5pt}\selectfont\bfseries\raggedright}{\theparagraph}{1em}{}
\titleformat{\subparagraph}[block]{\normalfont\fontsize{12pt}{14.5pt}\selectfont\bfseries\raggedright}{\thesubparagraph}{1em}{}
\titlespacing{\section}{0pt}{\baselineskip}{0.5\baselineskip}
\titlespacing{\subsection}{0pt}{\baselineskip}{0.5\baselineskip}
\titlespacing{\subsubsection}{0pt}{\baselineskip}{0.5\baselineskip}
\titlespacing{\paragraph}{0pt}{\baselineskip}{0.5\baselineskip}
\titlespacing{\subparagraph}{0pt}{\baselineskip}{0.5\baselineskip}
\numberwithin{equation}{section}

\newcommand{\footnoteONLYtext}[1]{%
 \begingroup
 \renewcommand{\thefootnote}{}%
 \footnotetext{#1}%
 \endgroup
}

\theoremstyle{plain}
\newtheorem{theorem}{Theorem}[section]
\numberwithin{theorem}{subsection}
\renewcommand{\thetheorem}{%
  \ifnum\value{subsection}>0 \thesubsection\else\thesection\fi.\arabic{theorem}}
\newtheorem{proposition}[theorem]{Proposition}
\newtheorem{lemma}[theorem]{Lemma}

\newtheorem{fact}[theorem]{Fact}
\theoremstyle{definition}
\newtheorem{definition}[theorem]{Definition}
\theoremstyle{remark}

\newcommand{\sphere}{\widehat{\mathbb C}}
\newcommand{\Crit}{\operatorname{Crit}}
\newcommand{\diam}{\operatorname{diam}}
\newcommand{\dist}{\operatorname{dist}}
\newcommand{\IR}{\operatorname{IR}}
\newcommand{\OR}{\operatorname{OR}}
\newcommand{\Phyp}{P_{\mathrm{hyp}}}
\newcommand{\Mhyp}{\mathscr M_{\mathrm{hyp}}}
\newcommand{\Pcl}{\mathcal P(f)}
\newcommand{\PC}{\operatorname{PC}(f)}
\newcommand{\tB}{\widetilde B}
\newcommand{\J}{J(f)}

\title{\normalfont\fontfamily{cmr}\selectfont\Large\bfseries Shrinking target problems for non-uniformly hyperbolic rational maps}
\author{Haoyang Ji}
\date{}

\hypersetup{
 pdftitle={Shrinking target problems for non-uniformly hyperbolic rational maps},
 pdfauthor={Haoyang Ji}
}

\begin{document}
\maketitle
\begin{abstract}% 
Let $f:\sphere\to\sphere$ be a rational map of degree at least two without
neutral cycles. For $y$ in the Julia set $J(f)$ and $\tau>0$, we study
the Hausdorff dimension of the shrinking target set
\[
 W_\tau(f,y):=\bigl\{x\in \sphere:d(f^n(x),y)<e^{-\tau n}
 \text{ for infinitely many }n\in\mathbb N\bigr\}.
\]
Assume that either (1) $f$ is an at most finitely renormalizable polynomial
with $|Df^n(f(c))|\to\infty$ for every critical point $c\in J(f)$,
or  (2) $\sum_{n\ge0}|Df^n(f(c))|^{-1}<\infty$ for every such $c$.
For every safe point $y\in J(f)$ and every $\tau>0$, we prove that
$\dim_H W_\tau(f,y)=s_\tau$, where $s_\tau$ is the unique solution
of the hyperbolic pressure equation $\Phyp(s)=\tau s$.
For the lower bound, we construct a Cantor subset of $W_\tau(f,y)$ using inverse branches associated with uniformly expanding sets. The key estimate gives a lower bound for the derivative of an inverse branch into the target ball, especially when pullbacks contain critical points.
For the Fibonacci quadratic polynomial $f$ and each $\tau>0$, we construct
an unsafe point $y\in J(f)\setminus\PC$ such that $\dim_H W_\tau(f,y)>s_\tau$.
\end{abstract}
\footnoteONLYtext{Date: \today}
\footnoteONLYtext{\textit{2020 Mathematics Subject Classification:}
37F35, 37C45, 37F10.}
\footnoteONLYtext{\textit{Keywords:} Shrinking targets, Hausdorff dimension, rational maps,
Julia sets, hyperbolic pressure.}

\section{Introduction}\label{sec:introduction}

Let $(X,d)$ be a separable metric space, $f:X\to X$ a Borel measurable
map and $\mu$ an ergodic $f$-invariant Borel probability measure on $X$.
Then {\it Birkhoff's ergodic theorem} implies that, for any
$y\in\operatorname{supp}\mu$,
\[
\liminf_{n\to\infty}d(f^n x,y)=0
\ \text{for $\mu$-almost every }x\in X.
\]
Thus, for $\mu$-typical $x\in X$, the orbit of $x$ visits
every neighborhood of $y$ infinitely many times. However, this statement
does not specify how fast the orbit approaches $y$.

To pass from qualitative to quantitative questions, let
$\psi:\mathbb{N}\to(0,\infty)$ satisfy $\psi(n)\to0$.
Here $\psi$ specifies the rate of approximation to the target.
Given $y\in X$, define the {\it shrinking target set}
\[
W(f,y,\psi):=\bigl\{x\in X:d(f^n x,y)<\psi(n)
\text{ for infinitely many }n\in\mathbb{N}\bigr\}.
\]
This set is a dynamically defined $\limsup$ set. Write
$B(y,r)=\{x\in X:d(x,y)<r\}$. For any $n\in\mathbb{N}$, let
\[
W_n(f,y,\psi):=f^{-n}(B(y,\psi(n))).
\]
Then the shrinking target set can be written as
\[
W(f,y,\psi)=\limsup_{n\to\infty}W_n(f,y,\psi)
=\bigcap_{N=1}^{\infty}\bigcup_{n\ge N}W_n(f,y,\psi).
\]
For $\tau>0$,  define the
{\it exponentially shrinking target set} by taking $\psi(n)=e^{-\tau n}$:
\[
 W_\tau(f,y):=\bigl\{x\in X:d(f^n(x),y)<e^{-\tau n}
 \text{ for infinitely many }n\in\mathbb N\bigr\}.
\]

Since $\mu$ is $f$-invariant, the $\mu$-measure of $W_n(f,y,\psi)$ is exactly the
$\mu$-measure of the ball $B(y,\psi(n))$. This together with the convergent
Borel--Cantelli lemma immediately implies that
\[
\sum_{n=1}^{\infty}\mu\bigl(B(y,\psi(n))\bigr)<\infty
\ \Longrightarrow\
\mu\bigl(W(f,y,\psi)\bigr)=0.
\]
The divergent case is considerably more delicate and typically requires quantitative
mixing or quasi-independence estimates. Under suitable assumptions, a dynamical
Borel--Cantelli theorem gives
\[
 \sum_{n=1}^{\infty}\mu(B(y,\psi(n)))=\infty
 \ \Longrightarrow\ \mu(W(f,y,\psi))=1.
\]
See, for example, Chernov and Kleinbock~\cite{CK01}
and Haydn {\it et al}.~\cite{HNPV13}.

The shrinking target problem has roots in metric Diophantine approximation.
For irrational $\alpha$, let $R_\alpha(x)=x+\alpha\pmod1$ on
$\mathbb R/\mathbb Z$ and let
$m$ denote normalized Lebesgue measure. Cassels~\cite{Cassels50}
proved that, for every $c>0$ and every $y\in\mathbb R/\mathbb Z$, $ m(W(R_\alpha,y,\psi_c))=1$ for $  \psi_c(n)=c/n$; see also Liao and Rams~\cite{LR13}.
Kurzweil~\cite{Kurzweil55}
proved that $\alpha$ is badly approximable if and only if $ m(W(R_\alpha,y,\psi))=1$ 
for every $y\in\mathbb R/\mathbb Z$ and every nonincreasing
$\psi$ as above with $\sum_{n=1}^{\infty}\psi(n)=\infty$;
see Simmons~\cite{Simmons} for a restatement.

When a shrinking target set is $\mu$-null, it is natural to ask for its
Hausdorff dimension. The dimension theory of shrinking target sets in
dynamical systems was initiated by Hill and Velani~\cite{HV}. For a rational
map $f$ that is uniformly expanding on its Julia set $J(f)$, they subsequently proved
in~\cite{HV97} that for every $y\in J(f)$ and $\tau>0$,
$\dim_H W_\tau(f,y)=s_\tau$, where $s_\tau$ is the unique solution of
\[
 P(f|_{J(f)},-s\log|Df|)=\tau s,
\]
and $P$ denotes topological pressure. This is a striking generalization
of Bowen's formula~\cite{Bowen79,Ruelle82}, which states that
$\dim_H J(f)$ is the unique solution of $P(f|_{J(f)},-s\log|Df|)=0$.

Stratmann, Urba\'nski and Zinsmeister~\cite[Theorem~1.1]{SUZ} considered
rational maps $f$ with parabolic periodic points for which
$\sum_{n\ge1}|Df^n(f(c))|^{-1}<\infty$ for every critical point $c\in J(f)$.
For a fixed parabolic periodic point $\omega$ and $\kappa>0$, they studied
the set $J_\omega^\kappa$ of points lying in infinitely many balls centered
at selected preimages of $\omega$. The ball centered at $z$ has radius
$(\rho_0|Df^n(z)|^{-1})^{1+\kappa}$, where $\rho_0>0$ is fixed and $z$ is
an $n$-th preimage of one of finitely many points in the backward orbit of
$\omega$. They proved that
\[
 \dim_H J_\omega^\kappa
 =\min\left\{\frac{D}{1+\kappa},
 \frac{D+\kappa p(\omega)}{1+\kappa(1+p(\omega))}\right\},
\]
where $D$ is the hyperbolic dimension of $f$ and $p(\omega)$ is the number
of attracting petals at $\omega$.

Dimension formulas for shrinking target sets are known for various dynamical
systems under suitable expansion assumptions.
The settings studied include countable expanding
Markov maps~\cite{Reeve11}, conformal iterated function systems~\cite{Urbanski02},
moving targets for continued fractions~\cite{LWWX14} and
$\beta$-transformations~\cite{BW14}, and matrix transformations of
tori~\cite{HV99,LLVZ23,LLVWZ25}. For exponentially shrinking moving
targets, Hu et al.~\cite{HPWZ25} obtained a dimension
lower bound for expanding toral endomorphisms and a dimension formula
for hyperbolic toral endomorphisms in dimension two.

Shrinking targets have also been studied in parameter space for families of
interval maps. For the real quadratic family $f_a(x)=1-ax^2$,
Bylund~\cite{Bylund23} proved that for Lebesgue a.e.
$a\in[1,2]$ for which $f_a$ has no attracting periodic orbit, $|f_a^n(0)|<1/n$ holds for infinitely many $n$. Aspenberg and
Persson~\cite{AP19} estimated the Hausdorff dimension of the set of
parameters $a$ for which the orbit of a prescribed point $X(a)$ within a family of 
piecewise expanding map $T_a$ hits shrinking targets around $y$ infinitely
often.

For further results, see the surveys
by He-Li-Velani~\cite{HLV25}, Hu-Huang-Li~\cite{HHL26} and
Athreya~\cite{Athreya09}.

In this article, we study the Hausdorff dimension of exponentially shrinking
target sets for non-uniformly hyperbolic rational maps on the Riemann sphere $\sphere$.
When $J(f)$ contains critical points, a target ball $B(y,e^{-\tau n})$
centered at $y\in J(f)$ may contain a critical value of $f^n$. In that
case, $f^n$ fails to be univalent on the pullback containing the
corresponding critical point. To control the effects of critical points,
we consider two conditions on derivative growth along the orbits of
critical values discussed in~\cite{RLS}. These are the
{\it large derivatives condition} for at most finitely renormalizable
polynomials and the {\it summability condition with exponent one} for
rational maps. Under either condition, the root of the hyperbolic
pressure equation gives a lower bound for $\dim_H W_\tau(f,y)$ for every
$y\in J(f)$ outside the postcritical set. Equality holds for every
{\it safe} target. For the Fibonacci quadratic polynomial $f$ and each
$\tau>0$, we construct an unsafe target $y\in J(f)\setminus\PC$ such that
$\dim_H W_\tau(f,y)>s_\tau$.

\subsection{Statement of results}\label{subsec:setting}

We identify the Riemann sphere $\sphere$ with
$\mathbb C\cup\{\infty\}$ and equip it with the spherical metric
whose density with respect to the Euclidean metric on $\mathbb C$
is $(1+|z|^2)^{-1}$. With this normalization,
$\diam(\sphere)=\pi/2$.
From now on, distances, balls, diameters, and derivative norms
are taken with respect to this metric.
For $z\in\sphere$ and $r>0$, let $B(z,r)$ denote the open ball
centered at $z$ with radius $r$.
We write $tV=B(x,tr)$ when $V=B(x,r)$ and $t>0$.
For nonempty subsets $A,B\subset\sphere$, write
$\dist(A,B)=\inf\{d(a,b):a\in A,\ b\in B\}$ and
$\dist(x,B)=\dist(\{x\},B)$. We use $a\asymp b$ to mean
$c_1b\le a\le c_2b$ for positive constants.

Let $\mathscr A$ denote the class of rational maps
$f:\sphere\to\sphere$ of degree at least two without neutral
periodic points. Let $\deg_z(f)$ denote the local degree of $f$ at $z$. We write
$\Crit(f)=\{z:\deg_z(f)>1\}$ for the set of critical points of $f$.

The {\it Julia set} $J(f)$ consists of the points $x\in\sphere$ for
which the family $\{f^n|_U\}_{n\ge0}$ is not equicontinuous on any
neighborhood $U$ of $x$. It is compact and completely invariant,
and is the closure of the repelling periodic points~\cite{Milnor}.

The {\it postcritical set} and its closure are defined as
\[
 \PC=\bigcup_{n\ge1}f^n(\Crit(f)), \ 
 \Pcl=\overline{\PC}.
\]

\begin{definition}\label{def:safe}
Following~\cite{P}, a point $y\in J(f)$ is called
{\it safe} if, for every $\sigma>0$, there is $n_0(\sigma)$ such that
\begin{equation}\label{eq:safe}
 \dist\left(y,\bigcup_{k=1}^n f^k(\Crit(f))\right)
       \ge e^{-\sigma n}, \ \forall n\ge n_0(\sigma).
\end{equation}
Every point of $J(f)\setminus\Pcl$ is safe, and every safe point lies
outside $\PC$.
\end{definition}

It follows directly from the definition that the set of points
in $J(f)$ which are not safe has Hausdorff dimension zero.

\begin{definition}
 A compact forward invariant set $K\subset J(f)$ is called
{\it uniformly expanding} if there are $C>0$ and $\lambda>1$ such that
$|Df^n(x)|\ge C\lambda^n$ for all $x\in K$ and $n\ge1$.
\end{definition}

A rational map is called {\it hyperbolic} or {\it expanding} if it is uniformly expanding
on its Julia set. For $f\in\mathscr A$, this holds
whenever $\Crit(f)\cap J(f)=\varnothing$.

Let $\mathscr H$ be the collection of all nonempty compact, forward
invariant and uniformly expanding subsets of $J(f)$.
Let $\mathscr H_{\mathrm{iso}}$ be the subcollection of isolated
elements of $\mathscr H$. Here $K$ is {\it isolated}
(or {\it forward locally maximal}) if there is an open neighborhood
$U$ of $K$ in $\sphere$ such that $f^n(x)\in U$ for all $n\ge0$ implies $x\in K$.
The elements of $\mathscr H_{\mathrm{iso}}$ are called
{\it expanding repellers}, following Ruelle.

For $s\ge0$,
define the {\it hyperbolic pressure} by
\[
 \Phyp(s)=\sup_{K\in\mathscr H}P(f|_K,-s\log|Df|),
\]
where $P(f|_K,\varphi)$ denotes the topological pressure of the
potential $\varphi$ on $K$. In~\cite[Definition~0.2]{PRLS04}, the authors define hyperbolic
pressure by
\[
 \widetilde P_{\mathrm{hyp}}(s)
 :=\sup_{K\in\mathscr H_{\mathrm{iso}}}
 P(f|_K,-s\log|Df|).
\]
By~\cite[Theorem~A]{PRLS04}, the two definitions give the same
pressure for $s\ge0$.

Let $\Mhyp$ be the collection of ergodic invariant probability
measures supported on members of $\mathscr H$. Write
$h_\mu=h_\mu(f)$ and $\chi_\mu=\int\log|Df|\,d\mu$ for the
{\it measure-theoretic entropy} and {\it Lyapunov exponent} of $\mu \in \Mhyp$, respectively. Then uniform expansion gives $\chi_\mu>0$ and the
variational principle yields
\begin{equation}\label{eq:variational}
 \Phyp(s)=\sup_{\mu\in\Mhyp}(h_\mu-s\chi_\mu).
\end{equation}
For $\tau>0$, put
\begin{equation}\label{eq:root-definition}
 s_\tau=\sup_{\mu\in\Mhyp}\frac{h_\mu}{\chi_\mu+\tau}.
\end{equation}

Our main result is the following.

\begin{theorem}\label{thm:main}
Let $f\in\mathscr A$ satisfy one of the following conditions:
\begin{enumerate}
\item $f$ is an at most finitely
renormalizable polynomial satisfying
\[
 \lim_{n \to \infty}|Df^n(f(c))| = \infty
 \ \text{for every }c\in\Crit(f)\cap J(f);
\]
\item $f$ satisfies
\[
 \sum_{n=0}^{\infty}\frac{1}{|Df^n(f(c))|}<\infty
 \ \text{for every }c\in\Crit(f)\cap J(f).
\]

\end{enumerate}
Then for every safe $y\in J(f)$ and every $\tau>0$,
\[
 \dim_H W_\tau(f,y)=s_\tau,
\]
where $s_\tau\in(0,2)$ is the unique solution of
\begin{equation}\label{eq:root}
 \Phyp(s)=\tau s.
\end{equation}
For every $y\in J(f)\setminus\PC$ and every $\tau>0$,
one still has
$\dim_H W_\tau(f,y)\ge s_\tau$.
\end{theorem}

Recall that a complex polynomial $f$ is
{\it renormalizable of period $p\ge1$} if there are Jordan disks
$U\Subset V$ in $\mathbb C$ for which $f^p:U\to V$ is proper
of degree at least two, the set
\[
 \{z\in U:f^{pn}(z)\in U\text{ for every }n\ge1\}
\]
is connected and properly contained in $J(f)$, and each critical
point belongs to at most one of $U,f(U),\ldots,f^{p-1}(U)$.
We call $f$ {\it at most finitely renormalizable} if only finitely
many integers $p\ge1$ satisfy these conditions.

Without the safety assumption, the lower bound can be strict,
as the following example shows.
Write $\omega(x)$ for the $\omega$-limit set of $x$.

\begin{theorem}\label{prop:unsafe}
Let $f$ be the Fibonacci quadratic polynomial defined in
Definition~\ref{def:fibonacci}. For every $\tau>0$,
there is a point $y\in\omega(0)\setminus\PC$ such that
\[
 \dim_H W_\tau(f,y)>s_\tau.
\]
The point $y$ is not safe and may depend on $\tau$.
\end{theorem}

% Planned subsection: Relation to previous results (shrinking targets,
% non-uniformly hyperbolic maps, dynamical Borel--Cantelli, hitting time statistics).

\subsection{Outline of the proof}

Let $f$ satisfy the hypotheses of Theorem~\ref{thm:main}.
Fix $y\in J(f)$ and $\tau>0$, put $A_n=f^{-n}(B(y,e^{-\tau n}))$. Then
\[
 W_\tau(f,y)=\limsup_{n\to\infty}A_n
 =\bigcap_{N=1}^{\infty}\bigcup_{n\ge N}A_n.
\]
In Proposition~\ref{prop:shrinking} we show that the diameters of
the components of $A_n$ tend uniformly to zero as $n\to\infty$.
Since each component intersects $J(f)$, we obtain $W_\tau(f,y)\subset J(f)$.

For the {\it upper bound}, suppose $y$ is safe. For all sufficiently large $n$,
the ball $B(y,2e^{-\tau n})$ contains no critical value of $f^n$.
Thus $A_n$ is the union of the univalent pullbacks $U_{n,z}$ of
$B(y,e^{-\tau n})$, where $U_{n,z}$ is the component of $A_n$
containing $z\in f^{-n}(y)$. The Koebe principle gives
$\diam U_{n,z}\le Ce^{-\tau n}|Df^n(z)|^{-1}$ with $C$ independent
of $n$ and $z$. The tree pressure formula~\cite[Theorem~2.2]{P} then gives
\[
 \limsup_{n\to\infty}\frac1n
 \log\sum_{z\in f^{-n}(y)}(\diam U_{n,z})^s
 \le \Phyp(s)-\tau s.
\]
For every $s>s_\tau$, the right-hand side is negative, so these covers
give $\mathcal H^s(W_\tau(f,y))=0$ and hence $\dim_H W_\tau(f,y)\le s_\tau$.

For the {\it lower bound}, fix $y\in J(f)\setminus\PC$ and $0<s<s_\tau$.
It suffices to construct a Cantor set $F\subset W_\tau(f,y)$ with $\dim_H F\ge s$.
The set $F$ can be realized as the limit set of a non-autonomous conformal iterated
function system on $V=B(y,a)$.

First, choose $\mu\in\Mhyp$ with entropy $h$ and Lyapunov exponent $\chi$
such that $h/(\chi+\tau)>s$.
Lemma~\ref{lem:free} applied to $\mu$ on its uniformly expanding support
then gives, for sufficiently small $\eta>0$, $M\ge2$ inverse branches
$F_i:2V\to V$ of a common iterate $f^q$ satisfying
\[
 \log M\ge q(h-\eta),\;
 e^{-q(\chi+\eta)}\le |DF_i(z)|\le\tfrac12
 \ \text{for }z\in\overline V,\ 1\le i\le M.
\]
Second, for a sufficiently small $\epsilon>0$, Proposition~\ref{prop:inverse}
gives, for every sufficiently small $r>0$,
an inverse branch $g_r$ of an iterate of $f$ on $2V$ satisfying
\[
 g_r(2V)\subset B(y,r),\;
 \inf_{\overline V}|Dg_r|\ge c_\epsilon r^{1+\epsilon},
\]
where $c_\epsilon>0$ is independent of $r$. The factor $r^\epsilon$
accounts for the loss in the derivative estimate when pullbacks contain critical points.

We shall choose integers $0=p_0<p_1<\cdots$ with large gap $p_j-p_{j-1}\ge2$.
At step $n$, the branch family of the IFS on $\overline V$ is
\[
 \mathcal F_n=
 \begin{cases}
 \{F_1,\ldots,F_M\}, & p_{j-1}<n<p_j,\\
 \{g_{r_j}\}, & n=p_j,
 \end{cases}
\]
where $r_j=e^{-\tau t_j}$, $B_j=B(y,r_j)$ and $t_j$ is the total iterate time (of $f$) of the first $p_j-1$ branches.
Each $p_j$ corresponds to a time $t_j$ at which every point of $F$ hits $B_j$, hence $F\subset W_\tau(f,y)$. The set $F$ supports a probability measure $\nu$ naturally, under which all cylinders of the same depth have equal mass.
By Lemma~\ref{lem:cylinders-to-balls} and the mass distribution principle,
it suffices to show that $\nu(Q)\le C_s(\diam Q)^s$ for every cylinder $Q$.

The only difficulty arises when the single branch $g_{r_j}$ is used, since
the cylinder keeps its mass while its diameter may decrease rapidly.
Let $Q$ be a cylinder at depth $p_{j-1}$ and
$P\subset Q$ a cylinder at depth $p_j$. With $k_j=p_j-p_{j-1}-1$, we have
$\nu(P)/\nu(Q)=M^{-k_j}\lesssim e^{-qk_j(h-\eta)}$.
By the chain rule and the Koebe principle, the derivative lower bounds give
\[
 \frac{\diam P}{\diam Q}
 \gtrsim e^{-qk_j(\chi+\eta)}r_j^{1+\epsilon}
 =e^{-qk_j(\chi+\eta)}e^{-(1+\epsilon)\tau t_j}.
\]
In Proposition~\ref{prop:target-criterion} we choose $p_j$ so that
$t_j/(qk_j)\to1$. Hence
\[
 \frac{\nu(P)}{(\diam P)^s}
 \lesssim
 e^{-qk_j[\,h-\eta-s(\chi+\eta+(1+\epsilon)\tau)\,]}
 \frac{\nu(Q)}{(\diam Q)^s}
 \ll \frac{\nu(Q)}{(\diam Q)^s},
\]
provided $\epsilon$ and $\eta$ are sufficiently small.
Lemma~\ref{lem:block-dimension} gives the corresponding estimate at every
depth and hence $\dim_H F\ge s$.

This paper is organized as follows.
In Section~\ref{sec:pullbacks} we use backward contraction to construct inverse
branches with arbitrary small images and prove the derivative estimate in
Proposition~\ref{prop:inverse}. Section~\ref{sec:free} combines these
branches with those from uniformly expanding sets to obtain the dimension
lower bound in Proposition~\ref{prop:target-criterion}. In
Subsection~\ref{subsec:safe-proof}, we use tree pressure to prove the upper
bound and complete the proof of Theorem~\ref{thm:main}. In
Subsection~\ref{subsec:unsafe}, we prove Theorem~\ref{prop:unsafe} by
constructing an unsafe target for the Fibonacci quadratic polynomial
for which the lower bound is strict.

\section{Backward contraction and inverse branches}\label{sec:pullbacks}

Throughout this section, let $f\in\mathscr A$. We first recall some background and definitions.

A {\it pullback} of a set $V$ by $f^n$, with $n\ge0$, is a
connected component of $f^{-n}(V)$, and $n$ is called the {\it pullback time}.

An {\it inverse branch} of $f^n$ on a domain $U$ is a
holomorphic map $\phi:U\to\sphere$ satisfying
$f^n\circ\phi=\mathrm{id}_U$. A preimage $z$ under $f^n$ is {\it regular}
if $Df^n(z)\ne0$.

The notion of {\it quasi-chain} was introduced in~\cite{BRLSvS08}.
We use the formulation from~\cite[Definition~4.2]{RLS}.

\begin{definition}\label{def:chains}
Let $f\in\mathscr A$, $n\ge0$, and let $\{G_k\}_{k=0}^n$ be a sequence
of nonempty connected open subsets of $\sphere$.
\begin{itemize}
\item[(1)] We call it a {\it chain} if $G_k$ is a connected component of
$f^{-1}(G_{k+1})$ for every $0\le k<n$.
\item[(2)] We call it a {\it quasi-chain} if
$G_k$ contains a connected component of $f^{-1}(G_{k+1})$
for every $0\le k<n$.
\end{itemize}
The {\it order} of a quasi-chain is
$\#\{0\le k<n:G_k\cap\Crit(f)\ne\varnothing\}$.
\end{definition}

\begin{definition}\label{def:expansion-away}
We say that $f$ is {\it expanding away from critical points} if
for every neighborhood $V$ of $\Crit(f)\cap J(f)$ that is open in
$\sphere$, the map $f$ is uniformly expanding on
\[
 K(V)=\{x\in J(f):f^n(x)\notin V\text{ for every }n\ge0\}.
\]
\end{definition}

Expansion away from critical points does not hold for every $f\in\mathscr A$. It does hold for polynomials that are at
most finitely renormalizable; see~\cite{KvS09}, and also~\cite{QY09}
for the totally disconnected case. 

We now recall the backward contraction property introduced
in~\cite{RL} for rational maps and in~\cite{BRLSvS08} for
interval maps. For $c\in\Crit(f)$ and $\delta>0$, we denote by $\tB(c,\delta)$ the connected component of
$f^{-1}(B(f(c),\delta))$ containing $c$.
For each $c\in\Crit(f)$, put $d_c=\deg_c(f)$. Since $\Crit(f)$ is finite, the local power form gives
constants $0<a_0\le a_1$ and $\delta_*>0$ depending only on $f$, such that for every $c\in\Crit(f)$ and
$0<\delta<\delta_*$,
\begin{equation}\label{eq:critical-neighborhood-bounds}
 B(c,a_0\delta^{1/d_c})\subset\tB(c,\delta)
 \subset B(c,a_1\delta^{1/d_c}).
\end{equation}

\begin{definition}\label{def:bc}
Let $R>1$. We say that $f$ is {\it backward contracting with
exponent $R$}, abbreviated $BC(R)$, if there is
$\delta_0=\delta_0(R)>0$ such that for every
$\delta\in(0,\delta_0)$, every $c\in\Crit(f)\cap J(f)$, every
$n\ge0$, and every pullback $W$ of $\tB(c,R\delta)$ by $f^n$,
\[
 \dist(W,f(\Crit(f)))\le\delta
 \ \Longrightarrow \ \diam W<\delta.
\]
The map $f$ is called {\it backward contracting} if it satisfies
$BC(R)$ for every $R>1$.
\end{definition}

Let $\mathscr A^{*}$ denote the class of maps
$f\in\mathscr A$ that are expanding away from
critical points and backward contracting.
The following fact gives
sufficient conditions for a map to belong to class $\mathscr A^{*}$.

\begin{fact}\label{fact:backward-contracting}
Let $f\in\mathscr A$ satisfy one of the following conditions:
\begin{enumerate}
\item $f$ is an at most finitely
renormalizable polynomial satisfying
\[
 \lim_{n \to \infty}|Df^n(f(c))| = \infty
 \ \text{for every }c\in\Crit(f)\cap J(f);
\]
\item $f$ satisfies
\[
 \sum_{n=0}^{\infty}\frac{1}{|Df^n(f(c))|}<\infty
 \ \text{for every }c\in\Crit(f)\cap J(f).
\]
\end{enumerate}
Then $f\in\mathscr A^{*}$.
\end{fact}

\begin{proof}
Backward contraction in \textnormal{(1)} and \textnormal{(2)}
follows from \cite[Theorem~A]{LS10} and
\cite[Theorem~A]{RL}, respectively.
Expansion away from critical points follows from
\cite{KvS09} or \cite[Corollary~8.3]{RL}
in \textnormal{(1)}, and from
\cite[proof of Lemma~3.1]{Prz98} in \textnormal{(2)}.
\end{proof}

The main result of this section is the following.

\begin{proposition}[Inverse branches with small images]\label{prop:inverse}
Let $f\in\mathscr A^{*}$.
Fix $y\in\J\setminus\PC$ and $\epsilon>0$.
There exists $a_*=a_*(f,y,\epsilon)\in(0,\pi/8)$ such that the following property holds for each $0<a<a_*$ and  $V=B(y,a)$. 

There are constants $0<r_0<a/2$ and $c_\epsilon,C_\epsilon>0$ depending only on
$f,y,\epsilon,a$ such that for
every $0<r<r_0$, there are an integer $N(r)\ge1$ and a univalent inverse
branch $g_r$ of $f^{N(r)}$ on $2V$ satisfying
\[
 g_r(2V)\subset B(y,r)\Subset V
\]
and
\begin{equation}\label{eq:weak-derivative}
 c_\epsilon r^{1+\epsilon}
 \le |Dg_r(z)|\le C_\epsilon r
 \ \text{for every }z\in\overline V.
\end{equation}

\end{proposition}

For hyperbolic rational maps, 
Hill and Velani in~\cite[Lemmas~3--4]{HV97} established inverse
branches $h_r$ for every sufficiently small $r>0$ satisfying
$cr\le |Dh_r(z)|\le Cr$ on disks of a fixed radius, where
$c,C>0$ are independent of $r$. In our setting, pullbacks may contain
critical points. Backward contraction controls the
losses in the derivative estimates, which are absorbed into the
factor $r^\epsilon$ in the lower bound of~\eqref{eq:weak-derivative}.

\subsection{Preparations}

We shall use the following spherical Koebe estimates.
The derivative distortion bound follows from~\cite[Lemma~1.2]{Prz98};
see also~\cite[Section~2.4]{PRL07}.
The ball inclusions and the distance estimate are standard geometric
consequences.

\begin{proposition}[Spherical Koebe Principle]\label{prop:spherical-koebe}
Let $f\in\mathscr A$. There are constants
$0<r_{\mathrm K}<\pi/4$ and $C\ge1$, depending only on $f$, such that for
every $y\in J(f)$, $0<r<r_{\mathrm K}$, $n\ge0$, and every inverse
branch $\phi$ of $f^n$ on $B(y,2r)$, we have:
\begin{equation}\label{eq:spherical-distortion}
 C^{-1}\le \frac{|D\phi(z)|}{|D\phi(w)|}\le C
 \ \text{for all }z,w\in B(y,r);
\end{equation}
for every $0<t\le r$,
\begin{equation}\label{eq:spherical-balls}
 B(\phi(y),C^{-1}t|D\phi(y)|)
 \subset \phi(B(y,t))
 \subset B(\phi(y),Ct|D\phi(y)|);
\end{equation}
and
\begin{equation}\label{eq:spherical-two-point}
 C^{-1}|D\phi(y)|  \le\frac{d(\phi(z),\phi(w))}{d(z,w)}\le C|D\phi(y)|
 \end{equation}
for all distinct $z,w\in\overline{B(y,r)}$.
\end{proposition}

We use the following consequence of {\it polynomial shrinking property} established in \cite{RLS}.

\begin{proposition}
\label{prop:shrinking}
Let $f\in\mathscr A^{*}$. There are constants
$0<\rho_{\mathrm{sh}}<\pi/2$ and $C_{\mathrm{sh}}>0$, depending only on $f$, such that
for every $z\in\J$, every $n\ge1$, and every pullback $U$ of
$B(z,\rho_{\mathrm{sh}})$ by $f^n$,
\begin{equation}\label{eq:shrinking}
 \diam U\le C_{\mathrm{sh}}n^{-1}.
\end{equation}
Moreover, for every $a>0$ there is a radius $0<\rho_a<\min\{\rho_{\mathrm{sh}},a/2\}$
such that, for every $z\in\J$, every $n\ge0$,
and every pullback $U$ of $B(z,\rho_a)$ under $f^n$, we have
$\diam U<a$.
\end{proposition}

\begin{proof}
The first assertion is exactly~\cite[Theorem~A]{RLS}
with polynomial shrinking exponent $1$.
For the second assertion, fix $a>0$ and choose an integer $N\ge1$ such that
$C_{\mathrm{sh}}/N<a$. Then every pullback $U$ of $B(z,r)$, with
$z\in\J$, $0<r\le\rho_{\mathrm{sh}}$, and $n\ge N$, satisfies $\diam U\le C_{\mathrm{sh}}/n\le C_{\mathrm{sh}}/N<a$.

It remains to treat the finitely many times $0\le j<N$.
For each such $j$, we claim that
\[
 \sup_{z\in\J}\;\sup_U\diam U\longrightarrow0
 \ \text{as }r\to0,
\]
where $U$ ranges over all pullbacks of $B(z,r)$ by $f^j$.

Otherwise, there exist $\epsilon_0>0$, a sequence $r_k\to0$,
points $z_k\in\J$, and pullbacks $U_k$ of $B(z_k,r_k)$ by $f^j$
such that $\diam U_k\ge\epsilon_0$ for every $k$.
Passing to a subsequence, we have $z_k\to z\in\J$, and the compact
connected sets $\overline{U_k}$ converge in the Hausdorff metric
to a compact connected set $F$ with $\diam F\ge\epsilon_0$.
By continuity of $f^j$, we have $f^j(F)=\{z\}$, so
$F\subset f^{-j}(z)$. Since $f^{-j}(z)$ is finite and $F$ is
connected, $F$ must be a singleton, contradicting
$\diam F\ge\epsilon_0$.

Since there are only finitely many times $0\le j<N$, we can choose
a common radius $0<\rho_a<\min\{\rho_{\mathrm{sh}},a/2\}$ such that all pullbacks
of $B(z,\rho_a)$ under $f^j$, $0\le j<N$, have diameter less than $a$,
uniformly for $z\in\J$. Together with the estimate for $n\ge N$,
this proves the second assertion.
\end{proof}

The following lemma shows that every pullback of a sufficiently
small spherical ball centered in $J(f)$ is simply connected.
A similar result was proved by Mihalache~\cite{Mihalache} for semi-hyperbolic rational maps.
The proof uses the following fact: if $D\subset\mathbb C$
is a simply connected domain,
and $P_d(z)=z^d$ for an integer $d\ge1$, then every connected component of
$P_d^{-1}(D)$ is simply connected.
This is a special case of~\cite[Proposition~3.1]{RGS}.
We apply this fact inductively, using the uniform control of
pullback diameters from Proposition~\ref{prop:shrinking}.

\begin{lemma}\label{lem:pullback-topology}
Let $f\in\mathscr A^{*}$.
There is a constant $0<r_{\mathrm{sc}}<\rho_{\mathrm{sh}}$,
depending only on $f$, such that for every $z\in\J$ and every
$0<r<r_{\mathrm{sc}}$, every pullback of $B(z,r)$ is simply connected.

\end{lemma}

\begin{proof}
For each $p\in\sphere$, choose topological disk neighborhoods $U_p$ of $p$
and $V_p$ of $f(p)$ such that $f:U_p\to V_p$ is conformally
equivalent to $w\mapsto w^{\deg_p(f)}$ on the unit disk.
By compactness, choose a finite subcover
$\{U_{p_1},\ldots,U_{p_m}\}$ of $\sphere$, and let $a_f>0$
be a Lebesgue number for this cover. Then every set
$E\subset\sphere$ with $\diam E<a_f$ is contained in
$U_{p_j}$ for some $1\le j\le m$.

Set $r_{\mathrm{sc}}=\rho_{a_f}$ as in Proposition~\ref{prop:shrinking}.
Then for every $z\in\J$ and $0<r<r_{\mathrm{sc}}$, every pullback of $B(z,r)$
has diameter less than $a_f$. Let $W$ be a pullback of $B=B(z,r)$ under $f^n$. Let $\{W_i \}_{i=0}^n$ be the chain with $W_0=W$ and $W_n=B$. Then $\diam W_i < a_f$ for $0 \leq i \leq n$. For each $i<n$, choose $U_p$ containing $\overline{W_i}$.
Then $W_{i+1}=f(W_i)\subset V_p$. Since $f:U_p\to V_p$ is conformally equivalent to a power map,
$W_i$ is simply connected whenever $W_{i+1}$ is. Starting
from the disk $W_n=B$, induction shows that $W=W_0$
is simply connected.

\end{proof}

The following lemma follows from~\cite[Lemma~4.2]{Bergweiler06}
and is a direct consequence of the Riemann--Hurwitz formula.

\begin{lemma}\label{lem:two-domains}
Let $F:U\to B$ be a proper holomorphic map of degree $D$ between
simply connected hyperbolic domains. Let $D_1,D_2\Subset B$ be Jordan
domains with disjoint closures.
Then there exist $i\in\{1,2\}$ and a component $W$ of
$F^{-1}(D_i)$ such that $F:W\to D_i$ is univalent.

\end{lemma}

\subsection{Controlled pullbacks and inverse branches}

Put $\ell=\max_{c\in\Crit(f)}\deg_c(f)$.
The following estimate is~\cite[Proposition~2]{LS} for $n\ge1$.
The case $n=0$ follows from~\eqref{eq:critical-neighborhood-bounds}.

\begin{lemma}\label{lem:critical-pullback-size}
Let $f\in\mathscr A^{*}$.
For every $0<\beta<1/\ell$, there are constants $C_\beta>0$ and
$t_\beta>0$, depending only on $f$ and $\beta$ such that for every
$c\in\Crit(f)\cap J(f)$, $0<t<t_\beta$, $n\ge0$, and every
pullback $Q$ of $\tB(c,t)$ by $f^n$, we have
\[
 \diam Q\le C_\beta t^\beta.
\]
\end{lemma}

For a domain $U\subset\sphere$ with
$\diam U<\pi/2$ and $x\in U$, define
\[
 \IR(U,x)=\inf_{z\in\partial U}d(x,z), \ 
 \OR(U,x)=\sup_{z\in\partial U}d(x,z).
\]
The following lemma follows from~\cite[Lemma~1.3]{Prz98}.
See also~\cite[Proposition~4]{LS} for assertion~\textnormal{(2)}
on small domains.

\begin{lemma}\label{lem:one-step-comparison}
Let $f\in\mathscr A$.
There is a constant $C_f\ge1$ depending only on $f$ such that
the following holds. Let $V\subset\sphere$ be a
domain and $U$ a component of $f^{-1}(V)$ with
$\diam U<\pi/2$ and $\diam V<\pi/2$.
\begin{itemize}
\item[(1)] We have
\[
 \sup_U|Df|\le C_f\frac{\diam V}{\diam U}.
\]
\item[(2)] For every $z\in U$,
\[
 \frac{\OR(U,z)}{\IR(U,z)}
 \le C_f\frac{\OR(V,f(z))}{\IR(V,f(z))}.
\]
\end{itemize}
\end{lemma}

\begin{proof}
By~\cite[Lemma~1.3]{Prz98}, there is $C_3\ge1$ depending only on $f$
such that 
\[
\sup_U|Df|\le C_3 \frac{\diam V}{\diam U}.
\]
It suffices to prove assertion (2). Fix $z\in U$ and choose $a\in\partial U$ with $d(z,a)=\IR(U,z)$.
Since $U$ is a component of $f^{-1}(V)$, we have $f(a)\in\partial V$.
The shortest spherical geodesic from $z$ to $a$ lies in $U$ except at its
endpoint. Integrating $|Df|$ along this geodesic gives
\[
 \IR(V,f(z))\le d(f(z),f(a))
 \le \sup_U|Df|\,\IR(U,z)
 \le C_3\frac{\diam V}{\diam U}\IR(U,z).
\]
Since $\diam V<\pi/2$, the antipode of $f(z)$ lies outside $\overline V$.
For any $u,v\in V$, we have
\[
 d(u,v)\le d(u,f(z))+d(f(z),v)\le2\OR(V,f(z)).
\]
Taking the supremum over $u,v\in V$ yields $\diam V\le2\OR(V,f(z))$.
Also, $\OR(U,z)\le\diam U$.
Hence
\[
 \frac{\OR(U,z)}{\IR(U,z)}
 \le 2C_3\frac{\OR(V,f(z))}{\IR(V,f(z))}.
\]
Taking $C_f=2C_3$ completes the proof.
\end{proof}

The following proposition first selects a pullback $W$ with
prescribed inner radius $\rho$. It also gives upper bounds for $\OR(W,x)$ and $\sup_W|Df^n|$.

\begin{proposition}\label{prop:controlled-pullbacks}
Let $f\in\mathscr A^{*}$ and put $b=1/(2\ell)$.
There are constants $\kappa_f\in(0,1]$, $L>1$ and $R_0>4^\ell$,
depending only on $f$, such that for every $R\ge R_0$ there is
a radius $\eta_R\in(0,\pi/4)$, depending only on
$f$ and $R$, with the following property.

For every $x\in\J$ and  every $0<\rho<\kappa_f\eta_R$, there exist
integers $n,m\ge0$ and $\eta'\in[\kappa_f\eta_R,\eta_R]$ such that
the component $W\ni x$ of $f^{-n}(B(f^n(x),\eta'))$ is simply
connected, has diameter less than $\pi/2$, and satisfies
\begin{equation}\label{eq:ir}
 \IR(W,x)=\rho,
\end{equation}
\begin{equation}\label{eq:order}
 m\le 1+\frac{\log(1/\rho)}{b\log(R/4^\ell)},
\end{equation}
\begin{equation}\label{eq:shape}
 \OR(W,x)\le L^{m+1}\rho,
\end{equation}
and
\begin{equation}\label{eq:forward-bound}
 \sup_W|Df^n|\le\frac{2\eta_R L^{m+1}}{\diam W}.
\end{equation}
Here $m$ is the order of the quasi-chain constructed in the proof.
\end{proposition}

Equation~\eqref{eq:ir} is essentially~\cite[Proposition~3]{LS}.
The proofs of~\eqref{eq:order}--\eqref{eq:forward-bound} follow
the quasi-chain argument used in~\cite[proof of Theorem~1]{LS}. We first fix some constants which depend only on $f$.

\begin{center}
\small
\renewcommand{\arraystretch}{1.25}
\setlength{\tabcolsep}{4pt}
\begin{tabular}{|p{2cm}|p{11cm}|}
\hline
Constant & Choice \\
\hline
$a_0,a_1$ & Local estimates in~\eqref{eq:critical-neighborhood-bounds} \\
\hline
$r_{\mathrm K},C$ & Koebe principle in Proposition~\ref{prop:spherical-koebe} \\
\hline
$r_{\mathrm{sc}}$ & Simple connectivity of pullbacks in Lemma~\ref{lem:pullback-topology} \\
\hline
$C_b,t_b$ & Pullback diameter bound in Lemma~\ref{lem:critical-pullback-size} with $\beta=b$ \\
\hline
$C_f$ & Derivative and shape estimates from Lemma~\ref{lem:one-step-comparison} \\
\hline
$h$ & $0<4h<\min\{r_{\mathrm K},r_{\mathrm{sc}}\}$ and
$h<\dist(c,J(f))$ for every critical point $c$ outside $J(f)$ \\
\hline
$M_f$ & $M_f>1$ with $a_0M_f^{1/d_c}\ge5a_1$ for every $c\in\Crit(f)$ \\
\hline
$\delta_*$ & Decrease $\delta_*$ so that $0<\delta_*<h$, $C_b\delta_*^b<1$, and (1)--(2) below hold \\
\hline
$\kappa_f$ & $(1+\max_{\sphere}|Df|)^{-1}$ \\
\hline
$K,L$ & $K=C^2$ and $L>KC_f$ \\
\hline
$R_0$ & $R_0=2\cdot4^\ell$ \\
\hline
\end{tabular}
\end{center}

For $0<t<\delta_*$, (1) follows from~\eqref{eq:critical-neighborhood-bounds}
and the choice of $M_f$, and (2) follows from Proposition~\ref{prop:shrinking}.
\begin{itemize}
\item[(1)] For every $c\in\Crit(f)$ and $z\in U\subset\tB(c,t)$,
\[
 B(z,2\diam U)\subset\tB(c,M_ft).
\]
\item[(2)] For every $c\in\Crit(f)\cap J(f)$, every pullback $Q$
of $\tB(c,t)$, including pullback time zero, satisfies $\diam Q<h$.
\end{itemize}
By Lemma~\ref{lem:pullback-topology}, every $Q$ in (2) is simply connected.
Every $Q$ in (2) contains a preimage of $c$ and hence meets $J(f)$.
By the choice of $h$, any critical point in $Q$ belongs to $J(f)$.

With the constants fixed above, the remaining choices and the construction proceed in the following order:
\[
\begin{gathered}
 \boxed{R\ge R_0}
 \ \longrightarrow\ \boxed{R_1,\ S,\ \delta_0(S)}
 \ \longrightarrow\ \boxed{\delta_R}
 \\[6pt]
 \longrightarrow\ \boxed{\zeta}
 \ \longrightarrow\ \boxed{\rho_{\zeta}}
 \ \longrightarrow\ \boxed{\eta_R}
 \\[6pt]
 \longrightarrow\ \boxed{n,\eta'}
 \ \longrightarrow\ \boxed{W}
 \ \longrightarrow\ \boxed{\{\widehat W_k\}_{k=0}^{n},\ m}
\end{gathered}
\]

\begin{proof}
For $R\ge R_0$, put
$R_1=R/4^\ell$ and $S=R\max\{1,M_f/4^\ell\}$.
Let $\delta_0(S)>0$ be given by $BC(S)$ and choose
\[
 0<\delta_R<\min\{\delta_0(S),t_b,\delta_*\}.
\]
Choose $0<\zeta<h$ such that $3\zeta<a_0\delta_R^{1/d_c}$ for every
$c\in\Crit(f)\cap J(f)$.
Take $\rho_\zeta$ from Proposition~\ref{prop:shrinking} and choose
$0<\eta_R<\rho_\zeta/2$.

Since $\rho_\zeta<\zeta/2<h<r_{\mathrm{sc}}$,
Proposition~\ref{prop:shrinking} and Lemma~\ref{lem:pullback-topology}
imply that for every $w\in J(f)$ and $0<r\le 2\eta_R$, every pullback $U$
of $B(w,r)$ is simply connected and has diameter less than $\zeta$.
Since $U$ contains a preimage of $w$ and $\diam U<\zeta<h$, the choice of $h$ gives
$U\cap\Crit(f)\subset J(f)$.
Whenever $c\in U\cap\Crit(f)$ and $z\in U$, we have
\begin{equation}\label{eq:first-critical-inclusion}
 B(z,2\diam U)\subset B(c,3\diam U)
 \subset B(c,a_0\delta_R^{1/d_c})\subset\tB(c,\delta_R).
\end{equation}
The constants $R_1,S,\delta_R,\zeta,\rho_{\zeta},\eta_R$ depend only on $f$ and $R$.

\noindent
{\bf Proof of \eqref{eq:ir}.} Given $x\in\J$ and
$0<\rho<\kappa_f\eta_R$, put
\[
 t_j=\max_{u\in\overline{B(x,\rho)}}d(f^j(u),f^j(x)), \ j\ge0.
\]
Then $t_0=\rho<\eta_R$. By the choice of $\kappa_f$,
$t_{j+1}\le\kappa_f^{-1}t_j$ for all $j \geq 0$.
We claim that there is an index $j$ with $t_j\ge\eta_R$. Otherwise, we have
$t_j < \eta_R<\rho_{\zeta} \leq \rho_{\mathrm{sh}}$ for all $j$. Hence the ball $B(x,\rho)$ would lie in a
pullback of $B(f^j(x),\rho_{\mathrm{sh}})$ for every $j\ge1$.
So~\eqref{eq:shrinking} gives
\[
 0<\diam B(x,\rho)\le C_{\mathrm{sh}}j^{-1},
\]
a contradiction as $j\to\infty$.
Let $n$ be smallest such that $t_{n+1} \geq \eta_R$ and set $\eta'=t_n$. Then
\[
 \kappa_f\eta_R\le\kappa_f t_{n+1}\le\eta'<\eta_R.
\]
By continuity and compactness, choose $\xi\in\overline{B(x,\rho)}$
such that
\[
 d(f^n(u),f^n(x))\le d(f^n(\xi),f^n(x))=t_n=\eta'
 \ \text{for every }u\in\overline{B(x,\rho)}.
\]
The open mapping theorem implies that $f^n(B(x,\rho))$ is open.
Thus $\xi\in\partial B(x,\rho)$ and
\[
 f^n(B(x,\rho))\subset B(f^n(x),\eta').
\]
Let $W$ be the component of $f^{-n}(B(f^n(x),\eta'))$ containing $x$.
Then $B(x,\rho)\subset W$ and $\xi\in\partial W$, so
$\IR(W,x)=\rho$.
Since $\eta'\le\eta_R$, the conclusion above shows that $W$ is simply connected and $\diam W<\zeta<\pi/2$.

\noindent
{\bf Proof of \eqref{eq:order}.}
For the rest of the proof, put $W_k=f^k(W)$ for $0\le k\le n$.
These domains form a chain with $\diam W_k<\zeta$
by the choice of $\eta_R$.
We define a quasi-chain $\{\widehat W_k\}_{k=0}^n$ by the following rules.
\begin{enumerate}
\item Set $\widehat W_n=B(f^n(x),2\eta_R)$.
\item For $k=n-1,\ldots,0$, let $\widehat W'_k$ be the component of
$f^{-1}(\widehat W_{k+1})$ containing $f^k(x)$, and set
\[
 \widehat W_k=
 \begin{cases}
 \widehat W'_k,
  &\text{if }\widehat W'_k\cap\Crit(f)=\varnothing,\\
 B(f^k(x),2\diam\widehat W'_k),
  &\text{otherwise}.
 \end{cases}
\]
\end{enumerate}
Inductively, $W_k\subset\widehat W'_k\subset\widehat W_k$ for $k<n$.
Let $m$ be the order of this quasi-chain. Put $k_0=n$, and let
$k_1>\cdots>k_m$ be all the indices in $\{0,\ldots,n-1\}$ for which
$\widehat W_{k_i}$ contains a critical point.

Assume $m\ge1$, and choose
$c_i\in\widehat W'_{k_i}\cap\Crit(f)$ for $1\le i\le m$.
The domain $\widehat W'_{k_1}$ is a pullback of
$B(f^n(x),2\eta_R)$, so $c_1\in J(f)$ by the choice of $\eta_R$.
Applying~\eqref{eq:first-critical-inclusion} with
$U=\widehat W'_{k_1}$ and $z=f^{k_1}(x)$ gives
\[
 \widehat W_{k_1}
 =B(f^{k_1}(x),2\diam\widehat W'_{k_1})
 \subset\tB(c_1,\delta_R).
\]
If $m\ge2$, then $\widehat W'_{k_2}$ is contained in a pullback
$Q$ of $\tB(c_1,\delta_R)$ by $f^{k_1-k_2}$.
Since $\delta_R<\delta_*$, (2) gives $\diam Q<h$. As $Q$ meets $J(f)$,
the choice of $h$ gives $c_2\in J(f)$.
The set $f(Q)$ is a pullback of $\tB(c_1,\delta_R)$ by
$f^{k_1-k_2-1}$ and contains $f(c_2)$.
Thus $BC(S)$ with $\delta=\delta_R/S<\delta_0(S)$ gives
$\diam f(Q)<\delta_R/S$.
Consequently,
$\widehat W'_{k_2}\subset\tB(c_2,\delta_R/S)$ and
\[
 \widehat W_{k_2}\subset\tB(c_2,M_f\delta_R/S)
 \subset\tB(c_2,\delta_R/R_1).
\]
Repeating the procedure, for each $1\le i\le m$,
\begin{equation}\label{eq:quasi-critical-scales}
 c_i\in J(f), \ 
 \widehat W_{k_i}\subset\tB(c_i,\delta_R R_1^{-(i-1)}).
\end{equation}
Since $f^{k_m}(W)\subset\tB(c_m,\delta_RR_1^{-(m-1)})$ and $\delta_RR_1^{-(m-1)}\le\delta_R<t_b$,
by Lemma~\ref{lem:critical-pullback-size} we have
\[
 \rho\le\diam W\le C_b\delta_R^b R_1^{-b(m-1)}\le R_1^{-b(m-1)},
\]
since $C_b\delta_R^b<C_b\delta_*^b<1$. Taking logarithms gives~\eqref{eq:order}.
For $m=0$, the same bound holds since $\rho<1$.

\noindent
{\bf Proof of \eqref{eq:shape}.}
Since $\delta_RR_1^{-(i-1)}\le\delta_R<\delta_*$,
condition \textnormal{(2)} in the choice of $\delta_*$ and
\eqref{eq:quasi-critical-scales} give
\[
 \diam\widehat W'_{k_i}
 \le \diam\tB(c_i,\delta_RR_1^{-(i-1)})<h,
\]
Thus each ball $\widehat W_{k_i}$ has radius
$2\diam\widehat W'_{k_i}<2h<r_{\mathrm K}$.
By construction, for $1\le i\le m$, the map
$F_i=f^{k_{i-1}-k_i-1}:W_{k_i+1}\to W_{k_{i-1}}$ extends to a
conformal map from $\widehat W_{k_i+1}$ onto $\widehat W_{k_{i-1}}$.
Similarly, $F_*=f^{k_m}:W\to W_{k_m}$ extends to a conformal map
from $\widehat W_0$ onto $\widehat W_{k_m}$.
For $m=0$, $f^n|_W$ extends to a conformal map from $\widehat W_0$
onto $\widehat W_n$.
These extensions are univalent because they have no critical points
and are proper maps onto simply connected balls.

We have $W_n\subset B(f^n(x),\eta_R)$ and, for $1\le i\le m$,
$W_{k_i}\subset\widehat W'_{k_i}\subset
B(f^{k_i}(x),\diam\widehat W'_{k_i})$.
By~\eqref{eq:spherical-balls} and~\eqref{eq:spherical-two-point},
for $m\ge1$
\[
 \frac{\OR(W_{k_i+1},f^{k_i+1}(x))}{\IR(W_{k_i+1},f^{k_i+1}(x))}
 \le K\frac{\OR(W_{k_{i-1}},f^{k_{i-1}}(x))}
                 {\IR(W_{k_{i-1}},f^{k_{i-1}}(x))}
 \ (1\le i\le m)
\]
and
\[
 \frac{\OR(W,x)}{\IR(W,x)}
 \le K\frac{\OR(W_{k_m},f^{k_m}(x))}{\IR(W_{k_m},f^{k_m}(x))}.
\]
For each $1\le i\le m$, Lemma~\ref{lem:one-step-comparison}\textnormal{(2)} gives
\[
 \frac{\OR(W_{k_i},f^{k_i}(x))}{\IR(W_{k_i},f^{k_i}(x))}
 \le C_f\frac{\OR(W_{k_i+1},f^{k_i+1}(x))}
                    {\IR(W_{k_i+1},f^{k_i+1}(x))}.
\]
Combining these estimates gives
\[
 \frac{\OR(W,x)}{\IR(W,x)}
 \le K^{m+1}C_f^m
       \frac{\OR(W_n,f^n(x))}{\IR(W_n,f^n(x))}
 =K^{m+1}C_f^m \leq L^{m+1}.
\]
This proves
\eqref{eq:shape}.

\noindent
{\bf Proof of \eqref{eq:forward-bound}.}
For $1\le i\le m$, let $H_i$ be the inverse of the conformal
extension of $F_i$ on $\widehat W_{k_{i-1}}$.
By~\eqref{eq:spherical-two-point},
\[
 \diam W_{k_i+1}
 \le C|DH_i(f^{k_{i-1}}(x))|\diam W_{k_{i-1}},
\]
while~\eqref{eq:spherical-distortion} gives
\begin{equation}\label{eq:block}
 \diam W_{k_i+1}\,\sup_{W_{k_i+1}}|DF_i|
 \le\frac{C\diam W_{k_i+1}}{|DH_i(f^{k_{i-1}}(x))|}
 \le K\diam W_{k_{i-1}}.
\end{equation}
The same argument for $F_*$ gives
\begin{equation}\label{eq:block-bottom}
 \diam W\,\sup_W|DF_*|\le K\diam W_{k_m}.
\end{equation}
For each $1\le i\le m$, Lemma~\ref{lem:one-step-comparison}\textnormal{(1)} gives
\begin{equation}\label{eq:single}
 \diam W_{k_i}\,\sup_{W_{k_i}}|Df|
 \le C_f\diam W_{k_i+1}.
\end{equation}
For $m\ge1$, the chain rule and~\eqref{eq:block}--\eqref{eq:single} give
\[
 \diam W\,\sup_W|Df^n|
 \le K^{m+1}C_f^m\diam W_{k_m}
 \prod_{i=1}^m\frac{\diam W_{k_{i-1}}}{\diam W_{k_i}}
 =K^{m+1}C_f^m\diam W_n.
\]
For $m=0$, we have $F_*=f^n|_W$ and $k_0=n$, so~\eqref{eq:block-bottom}
gives $\diam W\,\sup_W|Df^n|\le K\diam W_n$.
Since $\diam W_n\le2\eta_R$ and $K^{m+1}C_f^m\le L^{m+1}$,
we obtain~\eqref{eq:forward-bound}.
\end{proof}

\subsection{Proof of Proposition~\ref{prop:inverse}}

For fixed $f$ and $y$, the parameters in Proposition~\ref{prop:inverse}
are chosen in the following order, with $R$ chosen by~\eqref{eq:R-choice}.
\[
\begin{gathered}
 \boxed{\epsilon>0,\ a_\epsilon=1+\epsilon/2}
 \ \longrightarrow\ \boxed{R\ge R_0}
 \ \longrightarrow\ \boxed{\eta_R}
 \ \longrightarrow\ \boxed{\eta_*=\kappa_f\eta_R}
 \\[6pt]
 \longrightarrow\ \boxed{q_i,p_{i,\nu},n_{i,\nu},\kappa_{i,\nu}}
 \ \longrightarrow\ \boxed{a_*}
 \ \longrightarrow\ \boxed{0<a<a_*}
 \ \longrightarrow\ \boxed{m_\kappa,\ r_0}
 \\[6pt]
 \longrightarrow\ \boxed{0<r<r_0}
 \ \longrightarrow\ \boxed{\rho=r^{a_\epsilon}}
 \ \longrightarrow\ \boxed{n,\eta',W,m}
 \ \longrightarrow\ \boxed{N(r),g_r}
\end{gathered}
\]

\begin{proof}[Proof of Proposition~\ref{prop:inverse}]

Fix $y\in\J\setminus\PC$ and $\epsilon>0$.
Let $\kappa_f,L,R_0$ be given by 
Proposition~\ref{prop:controlled-pullbacks}. Put
$a_\epsilon=1+\epsilon/2$ and choose $R\ge R_0$ so large that
\begin{equation}\label{eq:R-choice}
 \frac{a_\epsilon\log L}{b\log(R/4^\ell)}<\epsilon/4.
\end{equation}
Let $\eta_R$ be given by
Proposition~\ref{prop:controlled-pullbacks} and put $\eta_*=\kappa_f\eta_R$.

Choose finitely many points $q_1,\ldots,q_M\in\J$ whose balls of
radius $\eta_*/4$ cover $\J$. In each $B(q_i,\eta_*/4)$ choose two distinct
points $p_{i,1},p_{i,2}$ and integers $n_{i,1},n_{i,2}\ge1$ such that
$f^{n_{i,\nu}}(p_{i,\nu})=y$ for $\nu=1,2$. Such points exist
because the backward orbit of $y$ is dense in $\J$ and $\J$ is perfect;
see~\cite{Milnor}. These preimages are regular because $y\notin\PC$.

For each $i$ and $\nu=1,2$,
let $\kappa_{i,\nu}$ be the local inverse branch of $f^{n_{i,\nu}}$
satisfying $\kappa_{i,\nu}(y)=p_{i,\nu}$.
Since there are finitely many branches and $p_{i,1}\ne p_{i,2}$ for each $i$,
choose $0<a_*<\pi/8$ small enough so that the following hold for every $0<a<a_*$.
Writing $V=B(y,a)$, each $\kappa_{i,\nu}$ is defined and univalent
on $4V=B(y,4a)$ and satisfies
\[
 f^{n_{i,\nu}}\circ\kappa_{i,\nu}=\mathrm{id}_{4V}, \ 
 \kappa_{i,\nu}(\overline{3V})\subset B(q_i,\eta_*/2)
\]
for $1\le i\le M$ and $\nu=1,2$. Moreover,
\[
 \kappa_{i,1}(\overline{3V})\cap\kappa_{i,2}(\overline{3V})
 =\varnothing \ \text{for every }1\le i\le M.
\]

Fix any $0<a<a_*$ and put
\[
 m_\kappa=\min_{i,\nu}\inf_{\overline V}|D\kappa_{i,\nu}|>0.
\]
This constant depends only on $f,y,\epsilon,a$.
Choose $r_0$ so that
\[
 0<r_0<\min\left\{ \frac{a}{2}, \eta_*^{1/a_\epsilon},
 L^{-8/\epsilon}\right\}.
\]
For each $0<r<r_0$, put $\rho=r^{a_\epsilon}<\kappa_f\eta_R$.
Apply Proposition~\ref{prop:controlled-pullbacks} with $x=y$
and this $\rho$ to obtain $n,\eta',W,m$.
By~\eqref{eq:order}, we have 
\[
 m+1\le 2+
 \frac{a_\epsilon\log(1/r)}{b\log(R/4^\ell)}.
\]
Since $0<r<1$,
by \eqref{eq:R-choice} we have
\[
 L^{m+1}
 \le L^2r^{-\frac{a_\epsilon\log L}{b\log(R/4^\ell)}}
 \le L^2r^{-\epsilon/4}
 \le r^{-\epsilon/2}.
\]
The pullback $W$ is simply connected. By~\eqref{eq:ir} and
\eqref{eq:shape}, we have $\diam W\ge\rho$ and
$\OR(W,y)\le L^{m+1}\rho\le r$. Thus~\eqref{eq:forward-bound} gives
\begin{equation}\label{eq:W-estimates}
 W\subset B(y,r), \ 
 \sup_W|Df^n|
 \le\frac{2\eta_R L^{m+1}}{\diam W}
 \le\frac{2\eta_R r^{-\epsilon/2}}{r^{1+\epsilon/2}}
 =2\eta_R r^{-1-\epsilon}.
\end{equation}
Put $W_n=f^n(W)=B(f^n(y),\eta')$, where
$\eta'\in[\eta_*,\eta_R]$.

Choose $i$ such that $d(f^n(y),q_i)<\eta_*/4$.
The Jordan domains $D_\nu=\kappa_{i,\nu}(3V)$, $\nu=1,2$,
have disjoint closures contained in $W_n$ (see Figure~\ref{fig:target-branches}).
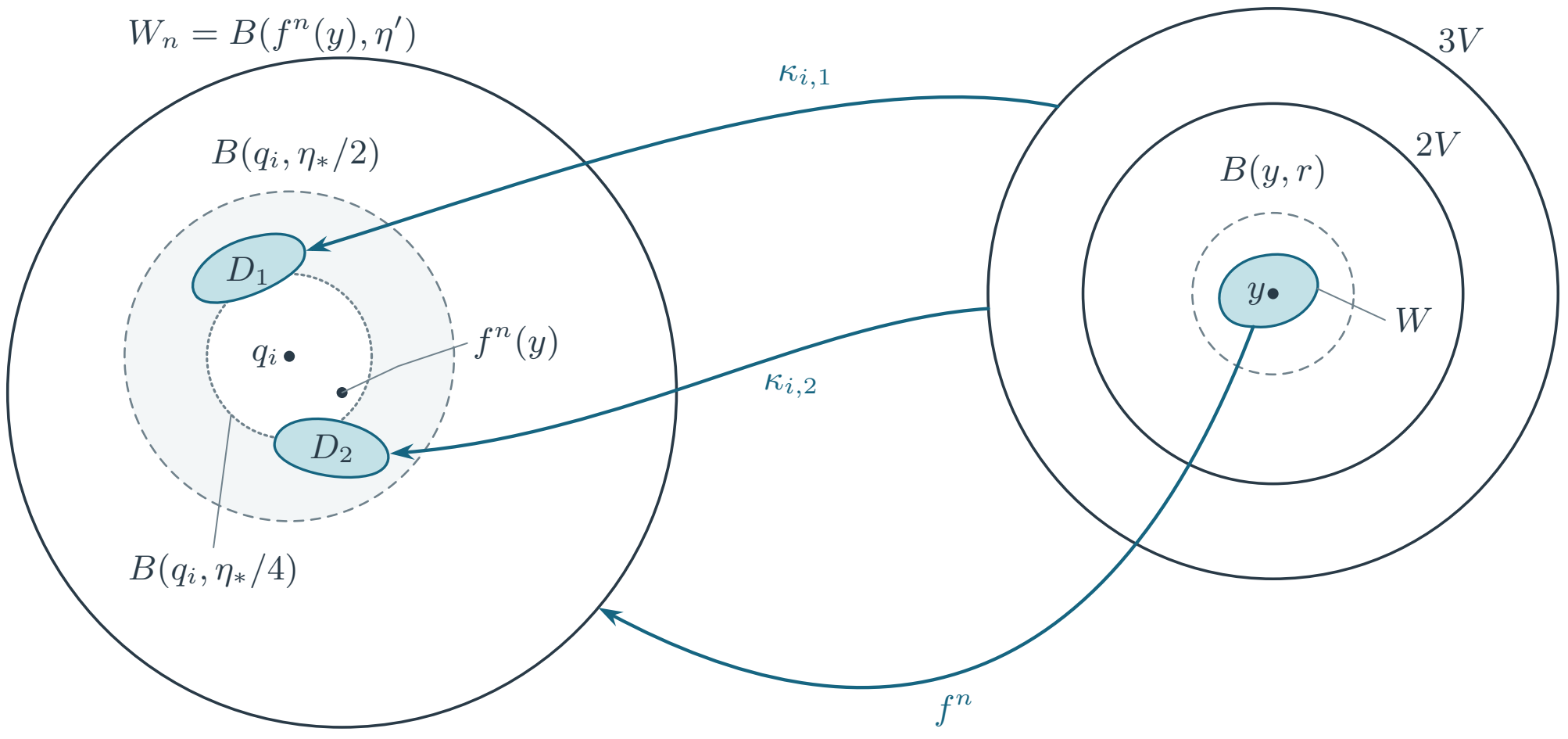
\begin{figure}[!htbp]
\centering
\begingroup
\setlength{\fboxsep}{7pt}
\resizebox{.85\linewidth}{!}{%
\colorbox{white}{%
\begin{tikzpicture}[
    scale=1.5,line cap=round,line join=round,font=\fontsize{9}{11}\selectfont,text=ink,
    boundary/.style={draw=ink,line width=.72pt},
    subdomain/.style={draw=guide,line width=.53pt,dash pattern=on 2.4pt off 2pt},
    pullback/.style={draw=accent,fill=imagefill,line width=.67pt},
    map/.style={draw=accent,text=accent,line width=.87pt,-{Stealth[length=2.2mm,width=1.45mm]}},
    every node/.style={inner sep=1.3pt}]
% Schematic radii. Named centres are shared by the circle and point commands.
\coordinate (fy) at (0,0);
\coordinate (q) at (-.32,.22);
\coordinate (y) at (5.65,.60);
\def\radiusWn{2.03}
\def\radiusQ{1.00}
\def\radiusThreeV{1.73}
\def\radiusTarget{.49}
\draw[boundary] (fy) circle (\radiusWn);
\draw[subdomain,fill=wash] (q) circle (\radiusQ);
\draw[draw=guide,densely dotted,line width=.65pt,fill=white]
  (q) circle[radius={.5*\radiusQ}];
\foreach \rad in {\radiusThreeV,2*\radiusThreeV/3}
  \draw[boundary] (y) circle[radius={\rad}];
\draw[subdomain] (y) circle (\radiusTarget);

% Three pullback domains, drawn independently of the metric circles.
\draw[pullback] (-.88,.58)
  .. controls (-.98,.69) and (-.78,.91) .. (-.51,.95)
  .. controls (-.33,.99) and (-.20,.93) .. (-.23,.82)
  .. controls (-.26,.68) and (-.76,.46) .. cycle;
\draw[pullback] (-.41,-.32)
  .. controls (-.41,-.17) and (-.18,-.12) .. (.04,-.19)
  .. controls (.29,-.26) and (.35,-.43) .. (.20,-.49)
  .. controls (.02,-.56) and (-.40,-.47) .. cycle;
\node at (-.57,.74) {$D_1$};
\node at (-.06,-.34) {$D_2$};

\begin{scope}[shift={(y)}]
  \draw[pullback] (-.30,-.12)
    .. controls (-.38,.01) and (-.27,.20) .. (-.06,.23)
    .. controls (.14,.27) and (.29,.16) .. (.27,.03)
    .. controls (.23,-.13) and (.04,-.22) .. (-.12,-.20)
    .. controls (-.21,-.20) and (-.27,-.17) .. cycle;
  \draw[guide,thin] (.27,.03) -- (.68,-.16)
    node[right=1pt,text=ink] {$W$};
\end{scope}

\foreach \point in {fy,q,y} \fill[ink] (\point) circle (.95pt);
\node[left=1pt] at (q) {$q_i$};
\draw[guide,thin] (fy) -- (.34,.16) -- (.76,.30)
  node[right,text=ink] {$f^n(y)$};
\node[left=.3pt] at (y) {$y$};

% Place labels in white space without breaking the boundary curves.
\node[above] at (-.42,2.03) {$W_n=B(f^n(y),\eta')$};
\node[above] at (-.28,1.30) {$B(q_i,\eta_*/2)$};
\node (innerlabel) at (-.78,-1.08) {$B(q_i,\eta_*/4)$};
\draw[guide,thin] ($(q)+(225:{.5*\radiusQ})$) -- (innerlabel.north);
\node[above right] at ($(y)+(56:\radiusThreeV)$) {$3V$};
\node[above right] at ($(y)+(44:{2*\radiusThreeV/3})$) {$2V$};
\node at ($(y)+(0,.74)$) {$B(y,r)$};

\draw[map] ($(y)+(139:\radiusThreeV)$)
  .. controls (3.04,2.00) and (1.18,1.29) .. (-.22,.86)
  node[pos=.35,above=4pt] {$\kappa_{i,1}$};
\draw[map] ($(y)+(183:\radiusThreeV)$)
  .. controls (2.86,.44) and (1.69,-.26) .. (.29,-.37)
  node[pos=.36,below=4pt] {$\kappa_{i,2}$};
\draw[map] ($(y)+(-.12,-.20)$)
  .. controls (4.62,-2.00) and (3.15,-2.16) .. (-40:\radiusWn)
  node[pos=.52,below=4pt] {$f^n$};
\end{tikzpicture}%
}%
}
\endgroup
\caption{The domains $D_1,D_2\subset W_n$ and the pullback $W\subset B(y,r)$.}
\label{fig:target-branches}
\end{figure}
By Lemma~\ref{lem:two-domains}, the proper map $f^n:W\to W_n$ has a
univalent inverse branch $h:D_\nu\to W$ for some $\nu\in\{1,2\}$.
Define $g_r=h\circ\kappa_{i,\nu}|_{2V}$, as
\[
 2V\xrightarrow{\kappa_{i,\nu}}D_\nu
    \xrightarrow{h}W\subset B(y,r).
\]
Set $N(r)=n+n_{i,\nu}$. Then
\[
 f^{N(r)}\circ g_r
 =f^{n_{i,\nu}}\circ\kappa_{i,\nu}|_{2V}
 =\mathrm{id}_{2V}.
\]

The chain rule and~\eqref{eq:W-estimates} give, for $z\in\overline V$,
\[
 |Dg_r(z)|
 =\frac{|D\kappa_{i,\nu}(z)|}{|Df^n(g_r(z))|}
 \ge\frac{m_\kappa}{2\eta_R}r^{1+\epsilon}.
\]
Take $c_\epsilon=m_\kappa/(2\eta_R)$.
Since $g_r(2V)\subset B(y,r)$, Cauchy's integral formula on $2V$
gives $|Dg_r|\le C_\epsilon r$ on $\overline V$ with $C_\epsilon$
depending only on $a$. Both bounds are uniform in $r,n$ and the
selected branches.
\end{proof}

\section{Hausdorff dimension of the constructed Cantor sets}\label{sec:free}

In this section, we construct Cantor subsets of $W_\tau(f,y)$
and estimate their Hausdorff dimensions from below.
The inverse branches obtained in Section~\ref{sec:pullbacks}
allow us to apply this construction to prove the lower bound
in Theorem~\ref{thm:main}.
In Subsection~\ref{subsec:branch-geometry}, we select a fixed family of
inverse branches from a uniformly expanding subset and prove a separation estimate. In
Subsection~\ref{subsec:cylinder-geometry}, we construct the Cantor sets
and estimate their dimensions using measures on cylinders.
Subsection~\ref{subsec:pressure-lower-bound} expresses the lower bound
in terms of hyperbolic pressure.

\subsection{Inverse branches}\label{subsec:branch-geometry}

We first state the consequences of~\cite[Lemma~11]{GPR} for rational
maps that are needed below.

\begin{lemma}\label{lem:gpr-branches}
Let $f\in\mathscr A$, and let
$\Lambda\subset J(f)$ be compact, invariant, uniformly expanding,
and topologically transitive. Let $\mu$ be an ergodic invariant
probability measure on $\Lambda$, with $h=h_\mu>0$ and
$\chi=\chi_\mu$, and let $O\subset\sphere$ be open with $\mu(O)>0$.
There is $\gamma_0>0$ such that for every $0<\gamma<\gamma_0$,
every $0<\epsilon<\min\{h,\chi\}$, and $\mu$-almost every
$z\in\Lambda$, there are $C\ge1$ and integers $n_k\to\infty$
with the following property.

For each $n=n_k$, there are $M_n$
inverse branches $\psi_{n,i}$ of $f^n$ on $U=B(z,\gamma)$ such that
$\psi_{n,i}(z)\in O\cap\Lambda$ and their images are pairwise disjoint.
For each $1\le i\le M_n$,
\[
 C^{-1}\le\frac{|D\psi_{n,i}(u)|}{|D\psi_{n,i}(v)|}\le C
 \ \text{for all }u,v\in U.
\]
Moreover, we have
\begin{equation}\label{eq:typical}
 \begin{gathered}
 M_n\ge C^{-1}e^{n(h-\epsilon)},\\
 C^{-1}e^{-n(\chi+\epsilon)}
 \le |D\psi_{n,i}(w)|\le Ce^{-n(\chi-\epsilon)}
 \ \bigl(w\in U,\ 1\le i\le M_n\bigr).
 \end{gathered}
\end{equation}
The constant $C$ may depend on $f,\Lambda,\mu,O,\gamma,\epsilon,z$,
but is independent of $n$ and $i$.
\end{lemma}

The estimate for $M_n$ follows from parts (iv) and (vii) of the
cited lemma; the remaining assertions follow from parts (i)--(iii)
and (v).

For each $K\in\mathscr H$, the variational principle gives
$h_{\mathrm{top}}(f|_K)=\sup_\mu h_\mu(f)$, where $\mu$ ranges over the
ergodic invariant probability measures on $K$.
By~\cite[Theorem~A]{PRLS04} at $s=0$,
$\sup_{K\in\mathscr H}h_{\mathrm{top}}(f|_K)=\log\deg f>0$.
Hence $\Mhyp$ contains a measure of positive entropy.

\begin{lemma}\label{lem:free}
Let $f\in\mathscr A$ and
$y\in J(f)\setminus\PC$. Fix $0<a_0<r_{\mathrm K}$, where $r_{\mathrm K}$
is given by Proposition~\ref{prop:spherical-koebe}, and put
$V_0=B(y,a_0)$ and $2V_0=B(y,2a_0)$.
Assume that, for every $r>0$, there are an integer $N = N(r)\ge1$
and an inverse branch $g = g_r$ of $f^N$ on $2V_0$ such that
\[
 g(2V_0)\subset B(y,r).
\]

Fix $\mu\in\Mhyp$ with $h=h_\mu>0$ and write $\chi=\chi_\mu$.
For every $0<\eta<\min\{h,\chi\}/4$, there are a radius $0<a<a_0$,
integers $q\ge1$, $M\ge2$, and $M$ inverse branches $F_1,\ldots,F_M$
of $f^q$ on $2V=B(y,2a)$ such that, with $V=B(y,a)$,
\begin{equation}\label{eq:free-images}
 \begin{gathered}
 F_i(2V)\subset V \ (1\le i\le M),\\
 F_i(\overline V)\cap F_j(\overline V)=\varnothing
 \ (1\le i<j\le M),
 \end{gathered}
\end{equation}
\begin{equation}\label{eq:free-estimates}
 \begin{gathered}
 \log M\ge q(h-\eta),\\
 e^{-q(\chi+\eta)}\le\inf_{\overline V}|DF_i|
 \le\sup_{\overline V}|DF_i|\le\tfrac12
 \ (1\le i\le M).
 \end{gathered}
\end{equation}
The radius $a$, the integers $q$ and $M$, and the inverse branches
$F_1,\ldots,F_M$ depend only on $f,y,a_0,\mu,\eta$.
\end{lemma}

\begin{proof}
Put $\Lambda=\operatorname{supp}\mu$. Since $\mu\in\Mhyp$, $\Lambda$ is
uniformly expanding. Ergodicity gives topological transitivity on $\Lambda$.
Moreover, $\mu$ is nonatomic since $h>0$.

{\it Step 1. Choose $\alpha$}. As $\PC$ is countable,
choose $z_0\in\Lambda\setminus\PC$. Since the backward orbit of $z_0$ is dense in $J(f)$, choose
$p\in B(y,a_0/2)$ and $t_1\ge1$ such that $f^{t_1}(p)=z_0$.
As $z_0\notin\PC$, we have $Df^{t_1}(p)\ne0$, so there exists an inverse branch $\alpha:U_0\to B(y,a_0/2)$ of $f^{t_1}$
on a neighborhood $U_0$ of $z_0$ with $\alpha(z_0)=p$.

Choose an open $O\Subset U_0$ with $\mu(O)>0$.
Fix $\epsilon=\eta/4$ and choose $0<\gamma<\gamma_0$ so that
$2\gamma<\dist(\overline O,\sphere\setminus U_0)$.
Since $\mu(O)>0$, choose $z\in O\cap\Lambda$ for which
Lemma~\ref{lem:gpr-branches} applies, and let $C$ and $\{n_k\}$
be the constant and sequence it gives. Put
\[
U=B(z,\gamma)\Subset U_0 \mbox{ and } E=\{w\in\sphere:\dist(w,\overline O)\le\gamma\}\subset U_0.
\]

{\it Step 2. Choose $\beta$}. Since $U$ meets $J(f)$ and the backward orbit of $y$ is dense in $J(f)$,
some preimage of $y$ lies in $U$. This preimage is regular because $y\notin\PC$.
Thus there are an integer $t_2\ge1$ and a sufficiently small
$0<a<a_0$ such that $V=B(y,a)\subset V_0$ and $f^{t_2}$ has
an inverse branch $\beta$ defined on a neighborhood of
$\overline{2V}$ with
\[
 \beta(\overline{2V})\subset U.
\]

{\it Step 3. Choose $g_*$}. Set $g_*=g_a$ and $t_*=N(a)$.
Then $g_*$ is an inverse branch of $f^{t_*}$ on $2V_0$ with
$g_*(2V_0)\subset V$.

With $\alpha$, $\beta$, and $g_*$ fixed, the derivatives of $\beta,\alpha,g_*$ are bounded above and away
from zero on $\overline V,E,\alpha(E)$, respectively.

{\it Step 4. Locate the images of $\psi_{n,i}$}. For each $n=n_k$,
let $\psi_{n,i}$, $1\le i\le M_n$, be the branches from
Lemma~\ref{lem:gpr-branches}. Integrating the derivative bound
in~\eqref{eq:typical} along radial arcs from $z$ in $U$ gives
\[
 \max_{1\le i\le M_n}\diam\psi_{n,i}(U)
 \le 2\gamma C e^{-n(\chi-\epsilon)}\longrightarrow0
 \ \text{as }n=n_k\to\infty.
\]
For all sufficiently large $k$ and every $i$, the image $\psi_{n,i}(U)$
has diameter less than $\gamma$. Since $\psi_{n,i}(z)\in O$,
we have $\psi_{n,i}(U)\subset B(\psi_{n,i}(z),\gamma)\subset E$.

Now we consider the following composition (see Figure~\ref{fig:free-branches}):
\[
 2V\xrightarrow{\beta}U\xrightarrow{\psi_{n,i}}U_0
    \xrightarrow{\alpha}2V_0\xrightarrow{g_*}V.
\]

\begin{figure}[!ht]
\centering
\begingroup
\setlength{\fboxsep}{7pt}
\resizebox{.85\linewidth}{!}{%
\colorbox{white}{%
\begin{tikzpicture}[
  font=\fontsize{9}{11}\selectfont, text=ink, line cap=round, line join=round,
  boundary/.style={draw=ink,line width=.72pt},
  subdomain/.style={draw=guide,line width=.53pt,dash pattern=on 2.4pt off 2pt},
  image/.style={draw=accent,fill=imagefill,line width=.67pt},
  branch/.style={draw=accent,line width=.87pt,-{Stealth[length=2.2mm,width=1.45mm]}},
  point dot/.style={circle,fill=ink,inner sep=0pt,minimum size=3pt},
  maplabel/.style={text=accent,inner sep=1pt}]

% Schematic radii. The four balls have the same centre.
\coordinate (y) at (0,0);
\draw[boundary] (y) circle (3.25);
\draw[subdomain] (y) circle (2.40);
\draw[boundary] (y) circle (1.65);
\draw[subdomain] (y) circle (.825);
\node at (0,2.81) {$2V_0=B(y,2a_0)$};
\node at (0,2.03) {$B(y,a_0/2)$};
\node at (0,1.24) {$2V$};
\node at (0,.34) {$V$};
\node[point dot] at (y) {};
\node[right=1pt] at (y) {$y$};
\draw[image,rotate around={12:(-.42,-.30)}] (-.42,-.30) ellipse (.205 and .11);
\coordinate (p) at (1.68,-1.12);
\draw[image,rotate around={-18:(p)}] (p) ellipse (.32 and .16);
\node[point dot] at (p) {};
\node at (.98,-2.62) {$p=\alpha(z_0)$};
\draw[guide,thin] (1.18,-2.43) -- (1.61,-1.32);

% Both z and psi_{n,i}(z) lie in O. The image domains are schematic.
\draw[boundary] plot[smooth cycle,tension=.85] coordinates
  {(4.72,-1.30) (4.87,1.30) (6.06,3.48) (8.46,4.00)
   (10.64,3.24) (11.49,1.23) (11.22,-1.07) (10.55,-3.26)
   (8.53,-3.92) (6.19,-3.55)};
\node at (10.00,2.64) {$U_0$};
\draw[subdomain,fill=wash] plot[smooth cycle,tension=.82] coordinates
  {(7.08,1.13) (7.98,1.68) (8.92,1.14) (9.03,.08)
   (8.82,-1.02) (8.18,-1.75) (7.55,-1.33) (7.21,-.30)};
\node at (8.75,-.30) {$O$};
\coordinate (z) at (7.53,1.13);
\draw[boundary] (z) circle (.88);
\node at (7.53,2.32) {$U=B(z,\gamma)$};
\node[point dot] at (z) {};
\node[right=1pt] at (z) {$z$};
\draw[image,rotate around={16:(7.18,.85)}] (7.18,.85) ellipse (.25 and .135);
\node[maplabel] at (6.30,1.79) {$\beta(2V)$};
\draw[guide,thin] (6.53,1.61) -- (7.12,1.02);
\coordinate (return) at (8.25,-1.51);
\draw[image,rotate around={-18:(return)}] (return) ellipse (.35 and .17);
\node[point dot] at (return) {};
\node at (9.20,-1.94) {$\psi_{n,i}(z)$};
\draw[guide,thin] (8.77,-1.82) -- (8.36,-1.58);
\coordinate (z0) at (6.52,-1.52);
\node[point dot] at (z0) {};
\node[above=2pt] at (z0) {$z_0$};

% The long arrows have separate upper and lower lanes.
\draw[branch] (1.62,.32) .. controls (3.58,.64) and (5.12,.91) .. (6.92,.85);
\node[maplabel] at (4.33,1.00) {$\beta$};
\draw[branch] (7.60,.25) .. controls (7.62,-.25) and (7.88,-.91) .. (8.18,-1.32);
\node[maplabel] at (8.11,-.34) {$\psi_{n,i}$};
\draw[branch] (6.44,-1.54) .. controls (5.18,-1.92) and (3.23,-1.62) .. (2.01,-1.20);
\node[maplabel] at (4.40,-2.02) {$\alpha$};
\draw[branch] (-3.15,-.80) .. controls (-2.49,-1.24) and (-1.42,-.97) .. (-.61,-.38);
\node[maplabel] at (-2.17,-1.49) {$g_a$};
\end{tikzpicture}%
}%
}
\endgroup
\caption{The inverse branches in $g_*\circ\alpha\circ\psi_{n,i}\circ\beta$.}
\label{fig:free-branches}
\end{figure}
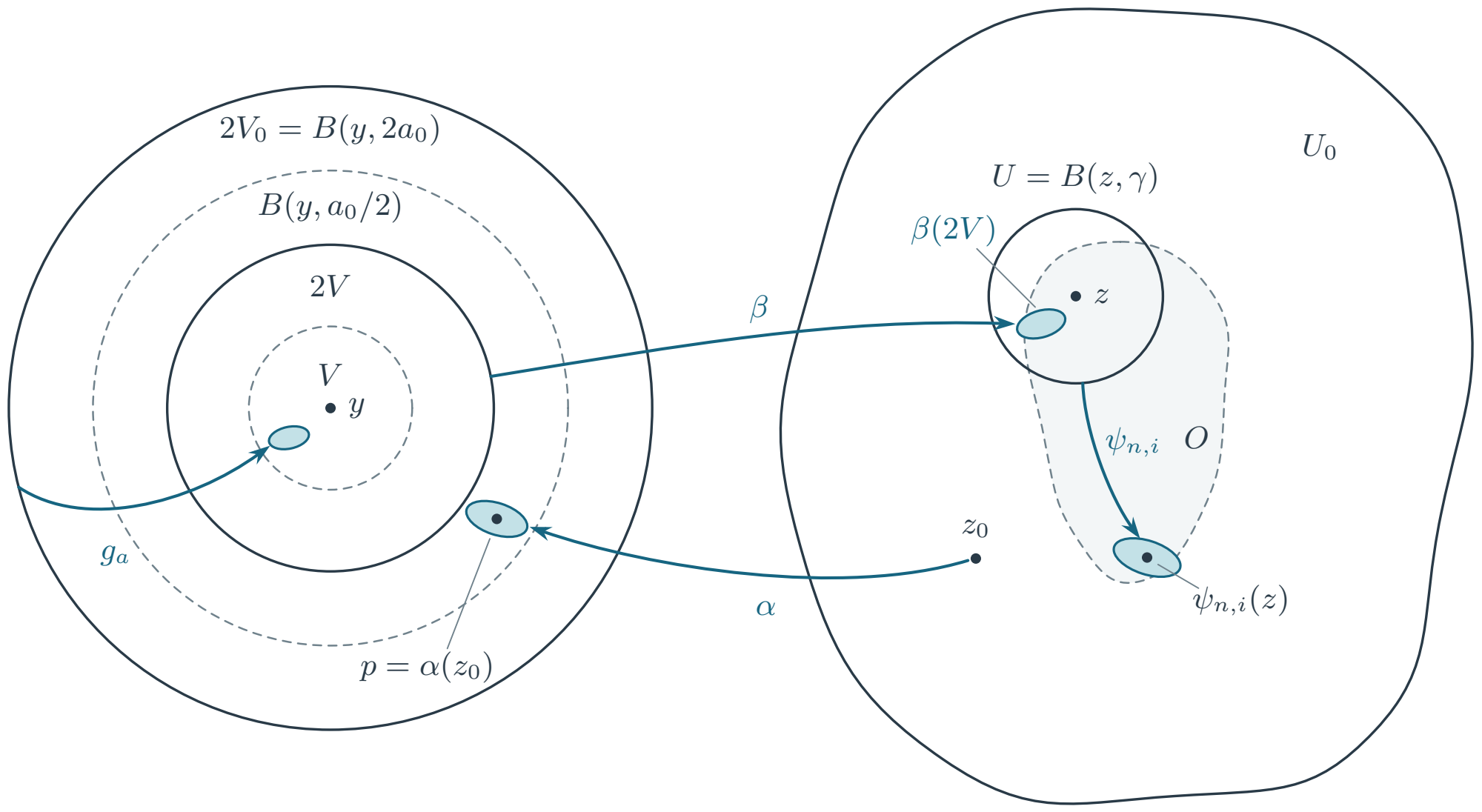
For each sufficiently large $n=n_k$, set $q_n=n+t_1+t_2+t_*$
and $F_{n,i}=g_*\circ\alpha\circ\psi_{n,i}\circ\beta$ for $1\le i\le M_n$.
Then $f^{q_n}\circ F_{n,i}=\mathrm{id}_{2V}$ and
$F_{n,i}(2V)\subset V$. For each such $n$, the images $F_{n,i}(\overline V)$
are pairwise disjoint, since the images of the $\psi_{n,i}$ are
disjoint and $g_*\circ\alpha$ is injective.

 Thus the chain
rule and~\eqref{eq:typical} give
\[
 C_1^{-1}e^{-n(\chi+\eta/4)}\le |DF_{n,i}|
 \le C_1e^{-n(\chi-\eta/4)}\ \text{on }\overline V.
\]
Here $C_1\ge1$ is independent of $n$ and $i$.
By~\eqref{eq:typical} and since $q_n-n=t_1+t_2+t_*$ is fixed,
\[
 \log M_n-q_n(h-\eta)
 \ge 3n\eta/4-\log C-(q_n-n)(h-\eta)\longrightarrow+\infty.
\]
Fix $n=n_k$ large enough that $\log M_n\ge q_n(h-\eta)>0$,
$\log C_1\le 3n\eta/4$, and $C_1e^{-n(\chi-\eta/4)}\le1/2$.
Then $M_n\ge2$ and for every $i$, we have
\[
 e^{-q_n(\chi+\eta)}\le e^{-n(\chi+\eta)}
 \le |DF_{n,i}|\le 1/2 \ \text{on }\overline V.
\]
Taking $q=q_n$, $M=M_n$, and $F_i=F_{n,i}$ gives~\eqref{eq:free-estimates}.

\end{proof}

\begin{lemma}\label{lem:geometry}
Let $f\in\mathscr A$ and $y\in J(f)\setminus\PC$.
Let $a$ and the inverse branches $F_1,\ldots,F_M$ be given by
Lemma~\ref{lem:free}, and put $V=B(y,a)$. There are constants
$c_0>0$ and $\delta>0$ such that every inverse branch $\Phi$
of an iterate of $f$ on $2V$, with $\Phi(2V)\subset V$, satisfies
\[
 \diam\Phi(\overline V)\ge c_0|D\Phi(y)|.
\]
Moreover, for $i\ne j$,
\begin{equation}\label{eq:relative-separation}
 \dist\bigl(\Phi(F_i(\overline V)),
             \Phi(F_j(\overline V))\bigr)
 \ge\delta\,\diam\Phi(\overline V).
\end{equation}
For fixed $f,V,F_1,\ldots,F_M$, the constants $c_0$ and $\delta$
are uniform in $\Phi$. 
\end{lemma}

\begin{proof}

Since $\diam\overline V=2a$, \eqref{eq:spherical-two-point} gives
\[
 (2a/C)|D\Phi(y)|\le\diam\Phi(\overline V)\le 2aC|D\Phi(y)|.
\]
Thus we may take $c_0=2a/C$.

Put 
\[
d_*:=\min_{i\ne j}\dist(F_i(\overline V),F_j(\overline V)) > 0.
\]
For $i\ne j$, the lower bound in~\eqref{eq:spherical-two-point} and the definition of $d_*$ give
\[
 \dist\bigl(\Phi(F_i(\overline V)),\Phi(F_j(\overline V))\bigr)
 \ge C^{-1}|D\Phi(y)|d_*
 \ge \frac{d_*}{2aC^2}\,\diam\Phi(\overline V).
\]
This proves~\eqref{eq:relative-separation} with $\delta=d_*/(2aC^2)$.
\end{proof}

\subsection{Cylinder geometry}\label{subsec:cylinder-geometry}

Suppose $f\in\mathscr A$ and $y\in J(f)\setminus\PC$. Fix $\mu\in\Mhyp$. Let $q,h,\chi,\eta$, $V=B(y,a)$ and $F_1,\ldots,F_M$ be as in Lemma~\ref{lem:free}.

For each $n\ge1$, take $\mathcal F_n=\{F_1,\ldots,F_M\}$ or
$\mathcal F_n=\{\phi\}$, where $\phi$ is an inverse branch of an iterate
of $f$ on $2V$ satisfying $\phi(2V)\subset V$ and
$\sup_{\overline V}|D\phi|\le1/2$.

Following~\cite[Definition~2.3]{RU16}, for each sequence
$\omega=(\phi_1,\ldots,\phi_k)$ with $\phi_j\in\mathcal F_j$, put
\[
 \Phi_\omega=\phi_1\circ\cdots\circ\phi_k, \ 
 Q_\omega=\Phi_\omega(\overline V).
\]
We call $Q_\omega$ a {\it cylinder of depth} $|\omega|=k$.
If each $\phi_j$ is an inverse branch of $f^{n_j}$, then $\Phi_\omega$
is an inverse branch of $f^{n_1+\cdots+n_k}$ and the {\it iterate time}
of $Q_\omega$ is $n_1+\cdots+n_k$.
The resulting {\it limit set} is
\[
 F=\bigcap_{k\ge1}
 \ \bigcup_{\omega\in\mathcal F_1\times\cdots\times\mathcal F_k}
 Q_\omega.
\]

For a cylinder of depth $k$, let $L=\#\{1\le n\le k:\#\mathcal F_n=M\}$.

\begin{lemma}\label{lem:cylinder-measure}
There is a Borel probability measure $\nu$ on $F$ such that for every
cylinder $Q_\omega$ of depth $k$,
\begin{equation}\label{eq:mass}
 \nu(Q_\omega)=\prod_{n=1}^{k}\frac{1}{\#\mathcal F_n}=M^{-L}.
\end{equation}
\end{lemma}

\begin{proof}
Let $\mathbb P_n$ be the uniform probability measure on $\mathcal F_n$.
On $\Sigma=\prod_{n\ge1}\mathcal F_n$, take the product measure
$\mathbb P=\otimes_{n\ge1}\mathbb P_n$.
Each map $\phi_j\in\mathcal F_j$ contracts distances on
$\overline V$ by at most $1/2$,
so $\diam Q_\omega\le2^{-|\omega|}\diam\overline V$.
Let $\pi(\{\phi_n\}_{n\ge1})$ be the unique
point in $\bigcap_{k\ge1}Q_{(\phi_1,\ldots,\phi_k)}$. By nesting and the uniform diameter upper bound, $\pi:\Sigma\to F$ is continuous.
Put $\nu=\pi_*\mathbb P$. By~\eqref{eq:free-images} and injectivity of the
inverse branches, distinct cylinders of the same depth are disjoint. Thus
$\pi^{-1}(Q_\omega)$ consists exactly of the sequences beginning with $\omega$.
By the definition of the product measure, this event has probability
$\prod_{n=1}^k(\#\mathcal F_n)^{-1}$, which gives~\eqref{eq:mass}.
\end{proof}

\begin{lemma}\label{lem:cylinders-to-balls}
Let $\delta$ be the separation constant from
Lemma~\ref{lem:geometry}, and let $\nu$ be a probability measure on the limit set $F$.
There exists an integer $N_\delta\ge1$, depending only on $\delta$,
such that the following holds.

Assume that for some $s>0$ and $C_s>0$, every cylinder $Q_\omega$ satisfies
\[
 \nu(Q_\omega)\le C_s(\diam Q_\omega)^s.
\]
Then
$\nu(B(x,r))\le N_\delta C_s r^s$ for every $x\in\sphere$ and
$0<r<\diam\overline V$. In particular, $\dim_HF\ge s$.
\end{lemma}

\begin{proof}
We use a stopping argument similar to that in~\cite[Theorem~3.2]{RU16}.
Fix $x\in\sphere$ and $0<r<\diam\overline V$.
For each infinite sequence $\{\phi_j\}_{j\ge1}$ with
$\phi_j\in\mathcal F_j$, let $k\ge1$ be the smallest integer
such that $\diam Q_{(\phi_1,\ldots,\phi_k)}\le r$, and select
$Q_{(\phi_1,\ldots,\phi_k)}$. Let $\mathcal S(r)$ be the
collection of all cylinders selected in this way. Note that $\mathcal S(r)$ is finite and covers $F$.

Let $Q,Q'\in\mathcal S(r)$ be distinct cylinders.
We show that $\dist(Q,Q')>\delta r$.
Write $Q=Q_{(\phi_1,\ldots,\phi_k)}$ and $Q'=Q_{(\psi_1,\ldots,\psi_l)}$.
Let $m$ be the largest integer with $0\le m < \min\{k,l\}$
such that $\phi_j=\psi_j$ for all $1\le j\le m$.
Put
\[
 P=Q_{(\phi_1,\ldots,\phi_m)},
\]
with $P=\overline V$ when $m=0$.
At step $m+1$, the two sequences choose distinct branches $\phi_{m+1}$ and $\psi_{m+1}$ from
$\{F_1,\ldots,F_M\}$. If $m=0$, then $Q$ and $Q'$ lie in distinct
first-level cylinders, whose distance is at least
$d_*\ge\delta\diam\overline V$ by the choice of $\delta$ in
Lemma~\ref{lem:geometry}. If $m\ge1$, put $\Phi=\phi_1\circ\cdots\circ\phi_m$,
so that $P=\Phi(\overline V)$. Then
\[
 Q  \subset \Phi\bigl(\phi_{m+1}(\overline V)\bigr), \ 
 Q' \subset \Phi\bigl(\psi_{m+1}(\overline V)\bigr).
 \]
Since $\diam P>r$ from the stopping rule, by~\eqref{eq:relative-separation}
\begin{equation}\label{eq:stop-separation}
 \dist(Q,Q')\ge  \dist( \Phi\bigl(\phi_{m+1}(\overline V)\bigr), \Phi\bigl(\psi_{m+1}(\overline V)\bigr) )  \geq \delta\diam P>\delta r.
\end{equation}

For each $Q\in\mathcal S(r)$ meeting $B(x,r)$, choose
$z_Q\in Q\cap B(x,r)$. By~\eqref{eq:stop-separation},
the balls $B(z_Q,\delta r/2)$ are pairwise disjoint
and lie in $B(x,(1+\delta/2)r)$. Spherical area comparison bounds
their number by an integer $N_\delta$ depending only on $\delta$.
By assumption we have
\[
 \nu(B(x,r))
 \le\sum_{\substack{Q\in\mathcal S(r)\\Q\cap B(x,r)\ne\varnothing}}
       C_s(\diam Q)^s
 \le N_\delta C_s r^s.
\]
The dimension lower bound follows from the mass distribution
principle~\cite[Chapter~4]{Falconer}.
\end{proof}

We now specify the steps \(n\) at which \(\mathcal F_n\) is a singleton and estimate the cylinder diameters to verify the assumption of Lemma~\ref{lem:cylinders-to-balls}.

Set $p_0=0$ and choose integers $p_j\ge p_{j-1}+2$ for $j\ge1$.
For each $j\ge1$, choose $g_j$ satisfying the conditions above for a singleton
branch and put
\begin{equation}\label{eq:cylinder-families}
 \mathcal F_n=
 \begin{cases}
 \{F_1,\ldots,F_M\}, & p_{j-1}<n<p_j,\\
 \{g_j\}, & n=p_j.
 \end{cases}
\end{equation}
Let $k_j=p_j-p_{j-1}-1$ and $K_j=\sum_{i=1}^j k_i$.
For a cylinder of depth $k$, take $j\ge0$ with $p_j\le k<p_{j+1}$.
Exactly $j$ of the first $k$ families are singletons, so $L=k-j =K_j+(k-p_j)$.

\begin{lemma}[A lower bound for the dimension]\label{lem:block-dimension}
Suppose the branches $g_j$ at steps $p_j$ satisfy, for some constants
$c,b>0$ and positive numbers $t_j$,
\[
 \inf_{\overline V}|Dg_j|\ge ce^{-bt_j} \ (j\ge1).
\]
Assume also that
\begin{equation}\label{eq:block-times}
 \sum_{i=1}^j t_i\le(1+o(1))qK_j, \ 
 j=o(K_j) \ (j\to\infty).
\end{equation}
Then
\[
 \dim_HF\ge\frac{h-\eta}{\chi+\eta+b}.
\]
\end{lemma}

\begin{proof}
Let $\nu$ be the measure from Lemma~\ref{lem:cylinder-measure} and fix a
cylinder $Q_\omega$ of depth $k$.
Applying Lemma~\ref{lem:geometry} to $\Phi_\omega$ and using the chain rule gives
\[
 \diam Q_\omega\ge c_0|D\Phi_\omega(y)|\ge c_0e^{-q(\chi+\eta)L}
                      c^j e^{-b\sum_{i=1}^j t_i}.
\]

Note that $L =K_j+(k-p_j) \geq K_j$. As $k\to\infty$, we have $j\to\infty$
and hence $L\to\infty$. By~\eqref{eq:block-times} we have
\[
\sum_{i=1}^j t_i\le qL+o(L) \mbox{ and } j=o(L) \mbox{ as } k \to \infty.
\] Taking logarithms in the
diameter estimate yields
\begin{equation}\label{eq:diam-final}
 -\log\diam Q_\omega\le q(\chi+\eta+b)L+o(L) \mbox{ as } k \to \infty.
\end{equation}

Fix $0<s<(h-\eta)/(\chi+\eta+b)$ and put
$\gamma=q(h-\eta-s(\chi+\eta+b))>0$.
Since $\log M\ge q(h-\eta)$, equations~\eqref{eq:mass}
and~\eqref{eq:diam-final} give
\[
 \log\frac{\nu(Q_\omega)}{(\diam Q_\omega)^s}
 \le-\gamma L+o(L) \mbox{ as } k \to \infty.
\]
For all sufficiently large $k$, every cylinder $Q_\omega$ of depth $k$
satisfies $\nu(Q_\omega)\le(\diam Q_\omega)^s$.
Only finitely many cylinders remain, so there is a constant
$C_s\ge1$ such that $\nu(Q_\omega)\le C_s(\diam Q_\omega)^s$ for every cylinder.
Lemma~\ref{lem:cylinders-to-balls} gives $\dim_HF\ge s$.
Letting $s\uparrow(h-\eta)/(\chi+\eta+b)$ proves the lemma.
\end{proof}

\subsection{Pressure and the dimension lower bound}\label{subsec:pressure-lower-bound}

For $f\in\mathscr A$, the function
$s\mapsto\Phyp(s)$ is nonincreasing and Lipschitz continuous;
see~\cite{PRLS04}. By~\eqref{eq:variational}, a Lipschitz
constant is $\log\max_{\sphere}|Df|>0$, since every
$\chi_\mu$, $\mu\in\Mhyp$, is positive and bounded by this constant.
Moreover, $\Phyp(0)>0$ because there are uniformly expanding sets with
positive entropy, while Ruelle's inequality
$h_\mu\le2\chi_\mu$ gives $\Phyp(2)\le0$;
see~\cite[Theorem~2]{Ruelle}.

For $\tau>0$, recall that
\begin{equation*}
 s_\tau=\sup_{\mu\in\Mhyp}\frac{h_\mu}{\chi_\mu+\tau}.
\end{equation*}

\begin{lemma}[Hyperbolic pressure roots]\label{lem:pressure-roots}
Let $f\in\mathscr A$.
For each $\tau>0$, the number $s_\tau$ belongs to $(0,2)$ and is the unique solution of $\Phyp(s)=\tau s$.
The function $\tau\mapsto s_\tau$ is continuous and strictly decreasing.
\end{lemma}

\begin{proof}
For $\tau>0$, the function $s\mapsto\Phyp(s)-\tau s$ is continuous,
strictly decreasing, and has opposite signs at $0$ and $2$, so it has a unique zero in $(0,2)$.
By~\eqref{eq:variational} and~\eqref{eq:root-definition},
\[
 \begin{aligned}
 \Phyp(s)-\tau s\le0
 &\iff h_\mu\le s(\chi_\mu+\tau) \ \text{for every }\mu\in\Mhyp\\
 &\iff s\ge\sup_{\mu\in\Mhyp}\frac{h_\mu}{\chi_\mu+\tau}=s_\tau.
 \end{aligned}
\]
Hence this zero is $s_\tau$.
For $0<\tau<\tau'$, we have
$\Phyp(s_\tau)-\tau's_\tau<0$, hence $s_{\tau'}<s_\tau$.
Since $\Phyp(s_{\tau'})\ge\Phyp(s_\tau)$, we have
\[
 0<s_\tau-s_{\tau'}
 \le\frac{s_{\tau'}}{\tau}(\tau'-\tau)
 \le\frac{2}{\tau}(\tau'-\tau).
\]
This proves continuity on $(0,\infty)$.
\end{proof}

The next proposition supplies the lower bound in
Subsection~\ref{subsec:safe-proof}; its construction is also used in
Subsection~\ref{subsec:unsafe}.

\begin{proposition}[Lower bound estimate]
\label{prop:target-criterion}
Let $f\in\mathscr A$ and
$y\in J(f)\setminus\PC$. Suppose there are a ball
$V_0=B(y,a_0)$ with $0<a_0<r_{\mathrm K}$, constants $c>0$, $\theta\ge1$ and $0<r_0<a_0$,
such that for every $0<r<r_0$ there are an integer $N(r)\ge1$ and an inverse
branch $g_r$ of $f^{N(r)}$ on $2V_0$ satisfying
\[
 g_r(2V_0)\subset B(y,r), \ 
 \inf_{\overline{V_0}}|Dg_r|\ge cr^\theta.
\]
Then $\dim_HW_\tau(f,y)\ge s_{\tau\theta}$ for every $\tau>0$.

\end{proposition}

\begin{proof}
Fix $\tau>0$ and $\mu\in\Mhyp$ with $h=h_\mu>0$ and $\chi=\chi_\mu$.
Since $g_r(2V_0)\subset B(y,r)$ for every $0<r<r_0$,
the inverse branch hypothesis of Lemma~\ref{lem:free} holds. Choose $0<\eta<\min\{h,\chi\}/4$
and apply that lemma to obtain $V=B(y,a)\subset V_0$ and
$q,M,F_1,\ldots,F_M$. Restrict the branches $g_r$ to $2V$. Cauchy's integral formula gives
$\sup_{\overline V}|Dg_r|\le Cr$, with $C$ independent of $r$.
Taking $r_0\le\min\{a,1/(2C)\}$, we have $g_r(2V)\subset V$ and
$\sup_{\overline V}|Dg_r|<1/2$ for every $0<r<r_0$.

We shall construct a limit set \(F \subset W_\tau(f,y)\) and choose integers \(1\le t_1<t_2<\cdots\) such that every \(z\in F\) satisfies
\[
f^{t_j}(z)\in B(y,e^{-\tau t_j}) \ \text{for every }j\ge1.
\]
For $j\ge1$, write $r_j=e^{-\tau t_j}$, $B_j=B(y,r_j)$ and $g_j=g_{r_j}$ whenever $r_j<r_0$.

{\it Construction of $F$.} Set $p_0=T_0=0$. Choose an integer $k_1\ge1$
large enough so that $e^{-\tau qk_1}<r_0$. Put
\[
p_1 = k_1 + 1, \ t_1 = q k_1, \ T_1 = t_1 + N(r_1).
\]
Inductively, for $j \ge 2$, choose an integer $k_j\ge1$ such that
\begin{equation}\label{eq:kj}
 qk_j\ge j^2\left(T_{j-1}+\sum_{i=1}^{j-1}t_i+j\right).
\end{equation}
Put
\[
 p_j=p_{j-1}+k_j+1, \ t_j=T_{j-1}+qk_j, \ T_j=t_j+N(r_j).
\]
It is clear that $\{t_j\}$ is increasing so that $r_j\le r_1<r_0$ and $N(r_j)$ is well-defined  before $k_{j+1}$ is chosen. Let $F$ be the limit set of the families in~\eqref{eq:cylinder-families},
\[
 F=\bigcap_{j\ge1}\ \bigcup_{|\omega|=p_j}Q_\omega.
\]
Each cylinder $Q_\omega$ contains $\Phi_\omega(y)\in J(f)$ by complete
invariance. Since the cylinder diameters tend to zero and $J(f)$ is
closed, $F\subset J(f)$.

{\it Hits at the times $t_j$.} For each $z\in F$, we prove
$f^{t_j}(z)\in B(y,e^{-\tau t_j})$ by induction. More precisely, we show that
$f^{t_j}\circ\Phi_{\omega_j}=g_j$ on $\overline V$ for every word
$\omega_j$ of length $p_j$. For $j=1$, write
$\Phi_{\omega_1}=H_1\circ g_1$, where $H_1$ is a composition of $k_1$
branches from $\{F_1,\ldots,F_M\}$. Since $t_1=qk_1$, we have
$f^{t_1}\circ H_1=\mathrm{id}$ on $V$, and hence
$f^{t_1}\circ\Phi_{\omega_1}=g_1$ on $\overline V$.

Suppose the identity holds at stage $j$. For each
$\omega_{j+1}=(\phi_1,\ldots,\phi_{p_{j+1}})$, put
$\omega_j=(\phi_1,\ldots,\phi_{p_j})$. Then
$\Phi_{\omega_{j+1}}=\Phi_{\omega_j}\circ H_{j+1}\circ g_{j+1}$,
where $H_{j+1}$ is the composition of the next $k_{j+1}$ branches.
By the induction hypothesis, $f^{T_j}\circ\Phi_{\omega_j}
=f^{N(r_j)}\circ g_j=\mathrm{id}$ on $\overline V$.
Since $t_{j+1}=T_j+qk_{j+1}$, we obtain
\[
 f^{t_{j+1}}\circ\Phi_{\omega_{j+1}}
 =f^{qk_{j+1}}\circ H_{j+1}\circ g_{j+1}
 =g_{j+1} \text{ on }\overline V.
\]
Thus every cylinder $Q_{\omega_j}$ of depth $p_j$ satisfies
$f^{t_j}(Q_{\omega_j})=g_j(\overline V)\subset B_j=B(y,e^{-\tau t_j})$.
These cylinders cover $F$ at each stage,
so every $z\in F$ hits the target at every $t_j$.

{\it Dimension estimate.} By~\eqref{eq:kj},
\[
 \sum_{i=1}^j t_i
 =qk_j+T_{j-1}+\sum_{i=1}^{j-1}t_i
 \le(1+j^{-2})qk_j
 \le(1+j^{-2})qK_j.
\]
Also $k_j\ge j^3/q$, so $j=o(K_j)$.
Moreover,
$\inf_{\overline V}|Dg_j|\ge cr_j^\theta
=ce^{-\theta\tau t_j}$.
Lemma~\ref{lem:block-dimension} therefore gives
\[
 \dim_HF\ge\frac{h-\eta}{\chi+\eta+\theta\tau}.
\]

Letting $\eta\downarrow0$ gives
$\dim_HW_\tau(f,y)\ge h_\mu/(\chi_\mu+\theta\tau)$.
Taking the supremum over $\Mhyp$ yields $s_{\tau\theta}$
by~\eqref{eq:root-definition}.
\end{proof}

\section{Proofs of Theorems~\ref{thm:main} and~\ref{prop:unsafe}}\label{sec:rational-proofs}

\subsection{Proof of the safe-target formula}\label{subsec:safe-proof}

For $f\in\mathscr A$, $n\ge1$, and $s>0$, define
\[
 Z_n(y,s)=\sum_{z\in f^{-n}(y)}|Df^n(z)|^{-s}
\]
for $y\in J(f)\setminus\bigcup_{j=1}^n f^j(\Crit(f))$. The sum runs over
the distinct preimages of $y$ in $\sphere$.
The following result relates these sums to hyperbolic pressure at safe points.

\begin{proposition}[Tree pressure {\cite[Theorem~2.2]{P}}]
\label{prop:tree-pressure}
Let $f\in\mathscr A$. For every safe
point $y\in J(f)$ and every $s>0$,
\[
 \limsup_{n\to\infty}\frac1n\log Z_n(y,s)=\Phyp(s).
\]
\end{proposition}

\begin{proof}[Proof of Theorem~\ref{thm:main}]
By Fact~\ref{fact:backward-contracting}, $f\in\mathscr A^{*}$.
Fix $\tau>0$. By Lemma~\ref{lem:pressure-roots}, $s_\tau\in(0,2)$
is the unique solution of the pressure equation~\eqref{eq:root}.

{\it The upper bound.}
Fix a safe point $y\in J(f)$ and $s>s_\tau$. Choose $\sigma\in(0,\tau)$
and set $a_n=\tfrac14e^{-\sigma n}$. Since
$\Phyp(s) <\tau s$, Proposition~\ref{prop:tree-pressure}
and~\eqref{eq:safe} allow us to take $N$ large enough that, for every $n\ge N$,
\[
 e^{-\tau n}<a_n<\min\{r_{\mathrm K},\rho_{\mathrm{sh}}\}, \ 
 Z_n(y,s)\le e^{n(\Phyp(s)+\tau s)/2},
\]
and $B(y,2a_n)$ contains no critical value of $f^n$.
For each $n\ge N$ and $z\in f^{-n}(y)$, let $\phi_{n,z}$ be the inverse branch
on this ball with $\phi_{n,z}(y)=z$, and put
$U_{n,z}=\phi_{n,z}(B(y,e^{-\tau n}))$.
Each $U_{n,z}$ lies in a pullback of $B(y,\rho_{\mathrm{sh}})$
by $f^n$, so~\eqref{eq:shrinking} gives
$\diam U_{n,z}\le C_{\mathrm{sh}}/n$, uniformly in $z$.
By Proposition~\ref{prop:spherical-koebe}, we have
\[
 \diam U_{n,z}\le C e^{-\tau n}|Df^n(z)|^{-1},
\]
where $C$ is independent of $n$ and $z$.
For each $n\ge N$, the sets $U_{n,z}$, $z\in f^{-n}(y)$, cover
$A_n=f^{-n}(B(y,e^{-\tau n}))$.
Since $W_\tau(f,y)\subset\bigcup_{n\ge N}A_n$, the family
$\{U_{n,z}:n\ge N,\ z\in f^{-n}(y)\}$ covers $W_\tau(f,y)$, and
\[
 \sum_{n\ge N}\sum_{z\in f^{-n}(y)}(\diam U_{n,z})^s
 \le C^s\sum_{n\ge N}e^{-\tau ns}Z_n(y,s)
 \le C^s\sum_{n\ge N}e^{-n(\tau s-\Phyp(s))/2}.
\]
This geometric sum tends to zero as $N\to\infty$, as does the diameter
bound $C_{\mathrm{sh}}/N$. Hence $\mathcal H^s(W_\tau(f,y))=0$
and $\dim_H W_\tau(f,y)\le s$. Letting $s\downarrow s_\tau$
gives the upper bound.

{\it The lower bound.}
For $y\in J(f)\setminus\PC$, apply
Proposition~\ref{prop:target-criterion} to the branches from
Proposition~\ref{prop:inverse}, with $\theta=1+\epsilon$, to obtain
$\dim_HW_\tau(f,y)\ge s_{\tau(1+\epsilon)}$ for every $\epsilon>0$.
Letting $\epsilon\downarrow0$ and using Lemma~\ref{lem:pressure-roots}
gives $\dim_HW_\tau(f,y)\ge s_\tau$.
Every safe point lies outside $\PC$, so the two bounds prove the theorem.
\end{proof}

\subsection{A counterexample for the Fibonacci quadratic polynomial}
\label{subsec:unsafe}

We first state the definition of Fibonacci combinatorics in quadratic polynomials.

\begin{definition}\label{def:fibonacci}
Let $q_1=1$, $q_2=2$, and $q_{j+1}=q_j+q_{j-1}$ for $j\ge2$.
The {\it Fibonacci quadratic polynomial} is the map
$f(z)=z^2+a_{\rm Fib}$, where $a_{\rm Fib}\in\mathbb R$ is the unique
parameter for which $f^4(0)<0$ and for each $j \geq 1$,
\[
 |f^{q_{j+1}}(0)|<|f^{q_j}(0)|.
\]
The integers $q_j$ are precisely the closest-return times of the
critical point $0$; see~\cite[Theorem~1.1]{LM}.
\end{definition}

Let $\Omega=\omega(0)$ be the $\omega$-limit set of the critical point $0$.
By~\cite{LM}, the critical point $0$ is recurrent and nonperiodic,
and $\Omega$ is a minimal Cantor set for $f$.
Moreover, $f\in\mathscr A$ and $0\in J(f)\setminus\PC$.

\begin{lemma}\label{lem:fibonacci-estimates}
Let $f(z)=z^2+a_{\rm Fib}$. For every $s>0$,
\[
 \sum_{n=1}^{\infty}|Df^n(f(0))|^{-s}<\infty.
\]
\end{lemma}
\begin{proof}
The estimate follows from~\cite[Lemma~5.9]{LM}, since along the bounded
orbit of $f(0)$ the spherical and Euclidean derivative norms of $f^n$
at $f(0)$ are comparable. See also
\cite[Corollary~1.1]{BLS} and~\cite{MN}.
Taking $s=1$ and applying Fact~\ref{fact:backward-contracting}(2) gives
$f\in\mathscr A^{*}$.
\end{proof}

\begin{lemma}\label{lem:fibonacci-safe}
The critical point $0$ of the Fibonacci quadratic polynomial is safe.
\end{lemma}
\begin{proof}
Put $d_j=d(0,f^{q_j}(0))$. Since spherical and Euclidean distances are
comparable near $0$, the estimate of Lyubich and Milnor~\cite[Theorem~1.3(1)
and Lemma~5.4]{LM} gives $-\log d_j=O(j^2)$.
As $q_j$ grows exponentially, $j^2=o(q_j)$. Fix $\sigma>0$.
For all sufficiently large $n$, take $j$ with $q_j\le n<q_{j+1}$.
Since the $q_j$ are the closest-return times of $0$,
\[
 \min_{1\le k\le n}d(0,f^k(0))=d_j\ge e^{-\sigma q_j}\ge e^{-\sigma n}.
\]
The only other critical point is the fixed point $\infty$, so~\eqref{eq:safe} holds.
\end{proof}

For the next two lemmas, fix $\tau>0$, $\tau/2<\alpha<\tau$,
and $\mu\in\Mhyp$ with $h_\mu>0$.
Write $h=h_\mu$ and $\chi=\chi_\mu$, and take $\epsilon>0$ and $0<\eta<\min\{h,\chi\}/4$.
Since $0\in J(f)\setminus\PC$,
we may fix $0<a_0<\min\{a_*(f,0,\epsilon),r_{\mathrm K}\}$ in Proposition~\ref{prop:inverse}
and put $V_0=B(0,a_0)$. For this $a_0$, the proposition gives
$r_0,c_\epsilon,C_\epsilon$ and inverse branches $g_r$ of $f^{N(r)}$
on $2V_0$ for $0<r<r_0$. Lemma~\ref{lem:free} gives a ball $V=B(0,a)\subset V_0$
and branches $F_1,\ldots,F_M$ of $f^q$ satisfying
\eqref{eq:free-images}--\eqref{eq:free-estimates}. Restrict $g_r$ to $2V$
and decrease $r_0$ so that $r_0<a$ and $C_\epsilon r_0<1/2$. Then
\[
 g_r(2V)\subset B(0,r)\subset V, \ 
 c_\epsilon r^{1+\epsilon}\le\inf_{\overline V}|Dg_r|
 \le\sup_{\overline V}|Dg_r|\le C_\epsilon r<1/2.
\]

As in the proof of Proposition~\ref{prop:target-criterion}, set $p_0=T_0=0$.
For a sequence of positive integers $\{k_j\}_{j\ge1}$, define inductively
\[
 p_j=p_{j-1}+k_j+1,\ 
 t_j=T_{j-1}+qk_j,\ 
 T_j=t_j+N(\tilde r_j),
\]
where $\tilde r_j=e^{-\alpha t_j}$ and each $k_j$ is chosen large enough that $\tilde r_j<r_0$.

\begin{lemma}\label{lem:unsafe-target}
There exist an unsafe point $y\in\Omega\setminus\PC$ and positive integers
$u_j<k_j$, $j\ge1$, such that \eqref{eq:block-times} holds, $\tilde r_j<r_0$ and
\begin{equation}\label{eq:fibonacci-target-image}
 f^{u_j}(B(0,\tilde r_j))\subset B(y,e^{-\tau(t_j+u_j)}).
\end{equation}
\end{lemma}

\begin{proof}
Since $f(z)-f(0)=z^2$, there are $A\ge1$ and $r_*>0$ such that
$d(f(z),f(0))\le A d(z,0)^2$ for $d(z,0)<r_*$ in the spherical metric.
With $L_f=\max\{2,\sup_{\sphere}|Df|\}$, we have
\[
 d(f^u(z),f^u(0))\le A L_f^u d(z,0)^2
 \ \text{for }u\ge1,\ d(z,0)<r_*.
\]

Choose a closed ball $D_0$ whose interior meets $\Omega$, and put $u_0=0$.
Since $\operatorname{int}D_{j-1}$ meets $\omega(0)$, the orbit of $0$ visits
$\operatorname{int}D_{j-1}$ infinitely often. Choose $u_j>\max\{u_{j-1},j\}$
with $v_j=f^{u_j}(0)\in\operatorname{int}D_{j-1}$. Since $0$ is recurrent
and nonperiodic, $v_j\notin\{f^i(0):1\le i\le j\}$.
Choose $k_j>u_j$ so large that $t_j=T_{j-1}+qk_j$ satisfies
\begin{equation}\label{eq:counter-growth}
 qk_j\ge j^2\Bigl(T_{j-1}+\sum_{i<j}t_i+u_j+j\Bigr),
\end{equation}
\[
\tilde r_j=e^{-\alpha t_j}<\min\{r_0,r_*\} \mbox{ and } A L_f^{u_j}e^{-2\alpha t_j}\le\tfrac14 e^{-\tau(t_j+u_j)}.
\]
The last inequality is possible because $2\alpha>\tau$.
Choose a closed ball $D_j\subset\operatorname{int}(D_{j-1})$ centered at
$v_j$ with radius at most $\tfrac14 e^{-\tau(t_j+u_j)}$ such that
$D_j\cap\{f^i(0):1\le i\le j\}=\varnothing$ and $\operatorname{int}D_j$ meets
$\Omega$. Set $T_j=t_j+N(\tilde r_j)$ before choosing $k_{j+1}$.

The balls $D_j$ shrink to a point $y\in\Omega$ since their centers belong
to $\Omega$. By construction, $y\ne f^i(0)$ for every $i\ge1$.
Since $\Omega\subset\mathbb R$ and $\PC=\{f^i(0):i\ge1\}\cup\{\infty\}$,
we have $y\notin\PC$.
By~\eqref{eq:counter-growth}, $t_j/u_j\ge j^2$. Hence
\[
d(y,f^{u_j}(0))\le\tfrac14 e^{-\tau(t_j+u_j)}
\le\tfrac14 e^{-\tau(j^2+1)u_j}<e^{-\sigma u_j}
\]
for every $\sigma>0$ and all sufficiently large $j$, which shows $y$ is unsafe. For every $z\in B(0,\tilde r_j)$,
\[
 d(f^{u_j}(z),y)\le A L_f^{u_j}\tilde r_j^2+d(v_j,y)
 \le\tfrac12 e^{-\tau(t_j+u_j)},
\]
giving~\eqref{eq:fibonacci-target-image}. Since~\eqref{eq:counter-growth}
implies~\eqref{eq:kj}, the calculation in the proof of
Proposition~\ref{prop:target-criterion} gives~\eqref{eq:block-times}.
\end{proof}

\begin{lemma}\label{lem:unsafe-construction}
For each $\epsilon>0$ and $0<\eta<\min\{h,\chi\}/4$ as above,
let $y$ be the unsafe point given by Lemma~\ref{lem:unsafe-target}.
Then there exists a compact Cantor set $F\subset W_\tau(f,y)$ such that
\[
 \dim_HF\ge\frac{h-\eta}{\chi+\eta+(1+\epsilon)\alpha}.
\]
\end{lemma}

\begin{proof}
Take $k_j$ as in Lemma~\ref{lem:unsafe-target} and put $g_j=g_{\tilde r_j}$ where $\tilde r_j=e^{-\alpha t_j}$.
Following the construction in the proof of Proposition~\ref{prop:target-criterion},
we obtain a Cantor set $F\subset W_\alpha(f,0)$. More precisely,
$f^{t_j}(F)\subset B(0,\tilde r_j)$ for every $j\ge1$.
By~\eqref{eq:fibonacci-target-image},
\[
 f^{t_j+u_j}(F)\subset f^{u_j}(B(0,\tilde r_j))\subset B(y,e^{-\tau(t_j+u_j)}).
\]
Since $t_j+u_j\to\infty$, we have $F\subset W_\tau(f,y)$.

Finally, $\inf_{\overline V}|Dg_j|\ge c_\epsilon \tilde r_j^{1+\epsilon}
=c_\epsilon e^{-(1+\epsilon)\alpha t_j}$. Lemma~\ref{lem:unsafe-target}
gives~\eqref{eq:block-times}, so Lemma~\ref{lem:block-dimension} with
$b=(1+\epsilon)\alpha$ gives the stated dimension bound.
\end{proof}

\begin{proof}[Proof of Theorem~\ref{prop:unsafe}]
Fix any $\tau>0$ and choose $\tau/2<\alpha<\tau$.
By Lemma~\ref{lem:pressure-roots}, $s_\alpha>s_\tau$,
so~\eqref{eq:root-definition} gives $\mu\in\Mhyp$ with
$h=h_\mu>0$, $\chi=\chi_\mu>0$, and $h/(\chi+\alpha)>s_\tau$.
Choose $\epsilon>0$ and $0<\eta<\min\{h,\chi\}/4$ small enough that
\[
 \frac{h-\eta}{\chi+\eta+(1+\epsilon)\alpha}>s_\tau.
\]
Take the unsafe target $y\in\Omega\setminus\PC$ from
Lemma~\ref{lem:unsafe-target}. For this $y$, Lemma~\ref{lem:unsafe-construction}
gives $\dim_H W_\tau(f,y)\ge\dim_H F>s_\tau$, proving the theorem.
\end{proof}

\section*{Funding}
This work was supported by the National Natural Science Foundation of China
[grant number 12301103].

\section*{Acknowledgements}
\noindent\textbf{AI declaration.}
The author used OpenAI Codex (GPT-6 Astra, GPT-6.1 Sol, and GPT-6.2 Astra)
to reorganize parts of the introduction and proofs, check how the constants
depend on the parameters, verify references, and turn the author's sketches into Latex and TikZ code for
Figures~\ref{fig:target-branches} and~\ref{fig:free-branches}.
The author directed and reviewed these revisions and takes responsibility
for the final text and figures.

\

\noindent\begin{minipage}{\textwidth}
\small
Haoyang Ji\\
School of Mathematics and Statistics, Zhengzhou University\\
Zhengzhou 450001, China\\
E-mail: \href{mailto:jihymath@zzu.edu.cn}{jihymath@zzu.edu.cn}
\end{minipage}
\end{document}